\documentclass[a4paper,11pt]{article}
\usepackage[english]{babel}
\usepackage[utf8]{inputenc}
\usepackage{german}
\usepackage{amsmath}
\usepackage{amsthm}
\usepackage{amssymb}
\usepackage{bbm}
\usepackage{graphicx}
\usepackage{pifont}
\usepackage{txfonts}
\usepackage{mathrsfs}
\usepackage{euscript}
\usepackage[ruled]{algorithm2e}
\usepackage{comment}
\usepackage{soul, color}
\usepackage{cmbright}
\usepackage[caption=false]{subfig}
\usepackage{float}
\usepackage[export]{adjustbox}
\usepackage[dvipsnames]{xcolor}
\usepackage[citecolor=NavyBlue,colorlinks=true,linkcolor=NavyBlue,urlcolor=NavyBlue,breaklinks]{hyperref}
\usepackage{url}

\DeclareMathAlphabet{\mathpzc}{OT1}{pzc}{m}{it}

\newcommand{\N}{\ensuremath{\mathbb{N}}}

\newcommand{\R}{\ensuremath{\mathbb{R}}}
\newcommand{\C}{\ensuremath{\mathbb{C}}}

\newcommand{\supp}{\mathrm{supp}}
\newcommand{\im}{\mathpzc{i}}
\newcommand{\dx}{\mathrm{d}}
\newcommand{\e}{\mathrm{e}}

\newcommand{\ds}{\displaystyle}
\newcommand{\ts}{\textstyle}

\newcommand{\re}{\mathrm{Re}}

\newcommand{\J}{\mathfrak{J}}

\theoremstyle{plain}
\newtheorem{theorem}{Theorem}[section]
\newtheorem{corollary}[theorem]{Corollary}
\newtheorem{lemma}[theorem]{Lemma}
\newtheorem{proposition}[theorem]{Proposition}
\newtheorem{definition}[theorem]{Definition}
\newtheorem{example}[theorem]{Example}
\newtheorem{remark}[theorem]{Remark}

\newcommand{\banf}{{\sf Proof.}\ }

\newcommand{\bend}{\hspace*{0ex} \hfill \hbox{\vrule height
    1.5ex\vbox{\hrule width 1.4ex \vskip 1.4ex\hrule  width 1.4ex}\vrule
    height 1.5ex}}

\newenvironment{Remark}{\goodbreak\begin{remark}\sl }{\end{remark}}

\newenvironment{Corollary}{\goodbreak\begin{corollary}\sl }{\end{corollary}}

\allowdisplaybreaks
\numberwithin{equation}{section}
\numberwithin{table}{section}
\numberwithin{figure}{section}
\title{Image Recovery for Blind Polychromatic Ptychography}
\author{Frank Filbir\thanks{Mathematical Imaging and Data Analysis, Helmholtz Center Munich, 85764 Neuherberg, Germany (\href{mailto:filbir@helmholtz-muenchen.de}{filbir@helmholtz-muenchen.de}, \href{mailto:oleh.melnyk@helmholtz-muenchen.de}{oleh.melnyk@helmholtz-muenchen.de}).}
\and Oleh Melnyk\footnotemark[1] \thanks{Department of Mathematics, Technical University of Munich, 85746 Garching bei M{\"u}nchen, Germany.} \thanks{Corresponding author.}}

\begin{document}
\maketitle
\begin{abstract}
Ptychography is a lensless imaging technique, which considers reconstruction from a set of far-field diffraction patterns obtained by illuminating small overlapping regions of the specimen. In many cases, a distribution of light inside the illuminated region is unknown and has to be estimated along with the object of interest. This problem is referred to as blind ptychography. While in ptychography the illumination is commonly assumed to have a point spectrum, in this paper we consider an alternative scenario with non-trivial light spectrum known as blind polychromatic ptychography. 

Firstly, we show that non-blind polychromatic ptychography can be seen as a recovery from quadratic measurements. Then, a reconstruction from such measurements can be performed by a variant of Amplitude Flow algorithm, which has guaranteed sublinear convergence to a critical point. Secondly, we address recovery from blind polychromatic ptychographic measurements by devising an alternating minimization version of Amplitude Flow and showing that it converges to a critical point at a sublinear rate. \\


\end{abstract}

\section{Introduction}\label{sec:Intro}
Ptychography is a scanning coherent diffraction imaging method and hence does not use advanced optical devices such as lenses for the image formation. 
The  image of the object is reconstructed numerically from data which consists of a stack of intensity measurements. This makes ptychographic imaging  predominantly a 
computational imaging technique.
The principle of ptychographic imaging can be outlined  as follows. An incoming coherent localized wave of a specific wavelength, in the
physics jargon called probe, is applied to illuminate a small region of the object of interest. The beam gets scattered and causes a diffraction pattern in the Fraunhofer or Fresnel region, depending on whether 
the diffraction plane is placed in the far- or near-field. A detector, usually a CCD camera, then records the intensity of that diffraction pattern. Subsequently the object 
is shifted to another position such that a different region will be illuminated by the localized beam and then the next measurement is recorded. In order to avoid a 
loss of information the adjacent illuminations have to overlap. 

\begin{figure}[h!]
    \centering
    \includegraphics
    [width = 1\linewidth ]
    {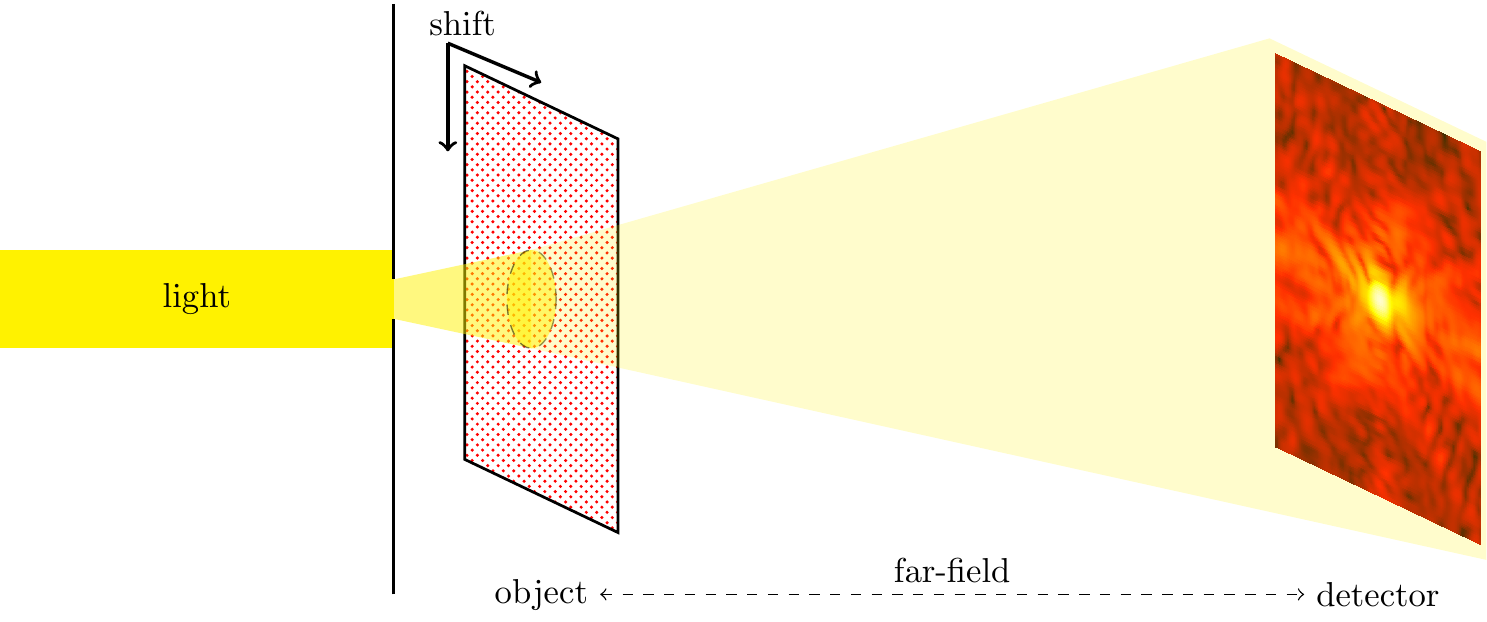}
    \caption{Ptychography.}
    \label{fig:1}
\end{figure}

In this way a set of intensity measurements is collected which forms the data base for the reconstruction process. The data redundancy allows to 
form an image of the object computationally. Over the last years the ptychographic technique was successfully used with different light sources such as 
synchrotron radiation \cite{Piazza.2014,Esmaeili.2015,Pfeiffer.2018}, electron beams \cite{Jiang.2018,Chen.2021}, and lasers \cite{Kharitonov.2021,Kharitonov.2022}. In a common experimental set-up light of one specific wavelength is used to illuminate the object, and the CCD camera is placed in the far-field (Fraunhofer distance). This experimental set-up then leads to measurements which are given 
mathematically as  
$$
\mathfrak{I}_{z}(x,\xi,\lambda)=\frac{1}{(\lambda z)^2}\left|\int_{\R^2}f(y)\, g_\lambda(y-x)\e^{-2\pi\im\frac{\xi\cdot y}{\lambda z}}\, \dx y\, \right|^2, 
$$ 
where $\lambda$ is the wavelength, $z$ is the distance of the object plane to the detector plane, $f$ is the object function, and $g_\lambda$ a wavelength dependent  window function 
which models the beam localization. In order to keep the exposition simple we will henceforth assume that $z$ is equal to one and the index $z$ will therefore be 
omitted. \\
The computational task now is the reconstruction of (an approximation) of $f$ from (samples) of $\mathfrak{I}(x,\xi,\lambda)$, i.e., from the squared absolute values of its Fourier transform. Hence the reconstruction problem is 
a phase retrieval problem. \\

As pointed out above, in the conventional experimental set-up one specific wavelength $\lambda$ is used, which then appears in the computational reconstruction 
process as a parameter. Moreover, often also the window function $g_\lambda$ is considered to be known beforehand. If, for example, the aperture has the form 
of a disc and the distance of the aperture to the object is sufficiently big, the window function is an Airy function which is frequently simply replaced by a 
Gaussian function. However, not every experimental configuration allows to have precise control over the window 
function. In those cases the window function has to be considered as an additional unknown object which we would like to  retrieve computationally as well. These 
category of problems are called blind ptychographic imaging and they were studied by several authors \cite{Maiden.2009,Thibault.2009,Hesse.2015,Chang.2019,Fannjiang.2020}. Giving up control about the concrete 
shape of the window function is however not the only necessary generalization of the problem. 
Light of only one specific wavelength is physically not easy to produce. Indeed, hard X-rays of one specific wavelength are usually produced by an 
electron-synchrotron, which is a huge machine. Other light producing systems however may provide light with a certain spectral distribution. Performing   
ptychographic measurements with spectrally distributed light will result in different intensity measurements. These are given in the form 
$$
\mathfrak{I}(x,\xi, \sigma)= \left| \int_\R \frac{1}{\lambda} \int_{\R^2}f_\lambda(y)\, g_\lambda(y-x)\e^{-2\pi\im\frac{\xi\cdot y}{\lambda}}\, \dx y \dx\sigma(\lambda) \, \right|^2, 
$$
where $\sigma$ is some compactly supported spectral density measure. Note that the object's scattering properties depends on the wavelength as well. Such model is considered in \cite{Burdet.2015} with an aim to improve quality of ptychographic reconstruction from light sources with near single wavelength illumination.  

If $\sigma$ consists of separated spectral lines represented by a weighed sum of Dirac's delta measures, i.e., $\sigma(\lambda) = \sum_{\ell = 1}^{L} \sigma_\ell \delta_{\lambda_\ell, \lambda}$ with $\lambda_\ell$ are enumerated such that $\lambda_1<\lambda_2<\dots<\lambda_L$, the intensity measurements reduce to
\begin{equation}\label{eq: measurements continuous}
\mathfrak{I}(x,\xi)=  \sum_{\ell=1}^{L} \frac{1}{\lambda_\ell^2} \big|\int_{\R^2}f_{\lambda_\ell}(y)\, w_{\lambda_\ell}(y-x) \e^{-2\pi\im\frac{\xi\cdot y}{\lambda_\ell}}\, \dx y\, \big|^2 , 
\end{equation}
where $w_{\lambda_\ell} = \sqrt{\sigma_\ell} g_{\lambda_\ell}$ is a window function for wavelength $\lambda_\ell$.

Recovery from noisy measurements of the form \eqref{eq: measurements continuous} is known as polychromatic ptychography. If the window functions $w_{\lambda_\ell}$ are unknown, similarly to the single wavelength case, the problem is referred to as blind polychromatic ptychographic imaging (BPPI). We note that the measurements \eqref{eq: measurements continuous} also arise if instead of polychromatic light multiple spatially separated apertures are used in ptychographic experiment \cite{Hirose.2020}. Moreover, a similar measurement model can be found in quantum state tomography \cite{Thibault.2013}.

In the literature, BPPI was addressed in several works \cite{Batey.2014, Guo.2018,Odstrcil.2018,Wei.2019,Metzler.2021}, each using a gradient-based method minimizing amplitude-based loss function. For instance, in \cite{Batey.2014} the authors establish a generalization of extended ptychographic iterative engine \cite{Maiden.2009} for polychromatic measurements \eqref{eq: measurements continuous} known as ptychographical information multiplexing method (PIM), which can be viewed as stochastic gradient descent applied to the amplitude-based loss function in analogy to \cite{Melnyk.2022}.
As the amplitude-based loss functions used for BPPI is non-Lipschitz, non-smooth and non-convex, convergence analysis of gradient-based methods for BPPI is not present in the literature. Furthermore, the absence of convergence guarantees leads to a non-trivial selection of the step sizes, which often requires multiple trial-and-error iterations to achieve good reconstruction.   

In this paper, we propose a new method for BPPI, which is based on alternating minimization technique \cite{Beck.2015}. In this way, we are able to reduce the reconstruction problem to repeated recovery from quadratic measurements \cite{Xu.2018b, Wang.2019, Huang.2020}. For such problems, it is possible to establish a gradient descent algorithm with appropriate step sizes, which guarantees sublinear convergence to a stationary point of the amplitude-based loss function similarly to \cite{Xu.2018}. Using these guarantees, we are also able to derive sublinear convergence of the whole alternating minimization technique.

The paper is structured as follows. Section \ref{sec:Pre} contains preliminaries about Wirtinger derivatives, gradient descent and its use for recovery from quadratic measurements. We return to the measurements \eqref{eq: measurements continuous} in Section \ref{sec:PP}. The established gradient optimization theory is applied first for non-blind problem and later for blind problem in Sections \ref{sec:NBP} and \ref{sec:BP}. In Section \ref{sec:NuEx} numerical trials are performed and proposed methods are compared to PIM. 


\section{Preliminaries}\label{sec:Pre}
In order to keep the presentation self-contained we start with some preliminary considerations regarding Wirtinger derivatives and a related gradient descent method. 
\subsection{Wirtinger derivatives and gradient descent}\label{sec:WDGD}
We start by collecting some facts on the Wirtinger calculus based on \cite{Hunger.2008,Bouboulis.2010}. Let $f(z)=u(x,y)+\im v(x,y),$ $z=x+\im y$ $x,y\in\R^n$ with real-valued differentiable functions $u$ 
and $v$. The function $f$ can be written as a function of the conjugate variables $z=x+\im y$ and $\bar{z}=x-\im y$. Since $u$ and $v$ are differentiable the 
function $f(z,\bar{z})$ is holomorphic w.r.t. $z$ for fixed $\bar{z}$ and vice versa. The Wirtinger calculus is a way to express
the derivatives of $f$ w.r.t. the real variables $x,y$ in terms of the conjugate variables $z$ and $\bar{z}$ treating them as independent.  
The so-called Wirtinger derivatives of $f$ are defined as   
\begin{equation}\label{eq:Pre0}
\partial_z f={\ts\frac{1}{2}}\,(\partial_x f -\im \partial_y f),\quad \partial_{\bar{z}} f={\ts\frac{1}{2}}\,(\partial_x f +\im \partial_y f).
\end{equation}
and we obviously have 
\begin{equation}\label{eq:Pre0a}
\overline{\partial_z f}=\partial_{\bar{z}}\bar{f},\quad\text{and}\quad \overline{\partial_{\bar z}f}=\partial_{z}\bar{f}.
\end{equation}
 The Wirtinger derivatives $\partial_zf$ and $\partial_{\bar{z}} f$ can also be expressed as 
\begin{align*}
\partial_z f
=\partial_z f(z,\bar{z})\big|_{\bar{z}=const.}=\begin{bmatrix} \partial_{z_1} f(z,\bar{z}),\dots, \partial_{z_n} f(z,\bar{z})\end{bmatrix}\big|_{\bar{z}=const.},\\[2ex]
\partial_{\bar{z}} f
=\partial_{\bar{z}} f(z,\bar{z})\big|_{z=const.}=\begin{bmatrix} \partial_{\bar{z}_1} f(z,\bar{z}),\dots, \partial_{\bar{z}_n} f(z,\bar{z})\end{bmatrix}
\big|_{z=const.}. 
\end{align*}
The Wirtinger gradient and Wirtinger Hessian are defined as 
\begin{equation*}
\nabla f(z)=\begin{bmatrix}(\partial_z f)^\ast\\ (\partial_{\bar{z}} f)^\ast\end{bmatrix},\quad 
\nabla^2f(z)=\begin{bmatrix}
                      \partial_{z}(\partial_z f)^\ast &\partial_{\bar{z}}(\partial_z f)^\ast\\[1ex]
                      \partial_z(\partial_{\bar{z}} f)^\ast&\partial_{\bar{z}}(\partial_{\bar{z}}f)^\ast
                     \end{bmatrix}.
\end{equation*} 
It follows immediately from \eqref{eq:Pre0a} that if $f$ is a real-valued function, i.e., $f(z)=u(x,y)$, the following relations hold  
\begin{equation}\label{eq:Pre1a}
\partial_{\bar{z}} f=\overline{\partial_z f},\quad \partial_{\bar{z}}(\partial_{\bar{z}}f)^\ast=\overline{\partial_z(\partial_z f)^\ast},\quad 
\partial_{z}(\partial_{\bar{z}} f)^\ast=\overline{\partial_{\bar{z}}(\partial_z f)^\ast}. 
\end{equation}
Henceforth we will use the less clumsy notation 
$$
\nabla_z f:=(\partial_z f)^\ast,\quad \nabla^2_{z,z} f:=\partial_{z}(\partial_z f)^\ast,\  \text{resp.} \quad \nabla^2_{z,\bar{z}} f:=\partial_{z}(\partial_{\bar{z}} f)^\ast.
$$
The second-order Taylor polynomial of $f$ at a point $z_0$ reads as 
$$
P_f(u,z_0)=f(z_0)+(\nabla f(z_0))^\ast \begin{bmatrix} u\\ \bar{u}\end{bmatrix}\, +\, \begin{bmatrix} u\\ \bar{u}\end{bmatrix}^\ast\, \nabla^2f(z)\, 
\begin{bmatrix} u\\ \bar{u}\end{bmatrix}
$$
In case of a real-valued function $f$, the quadratic term of the Taylor polynomial $P_f$ can be expressed in the following way 
\begin{equation}\label{eq:Pre2 real val}
\begin{bmatrix} u\\ \bar{u}\end{bmatrix}^\ast\, \nabla^2f(z)\, \begin{bmatrix} u\\ \bar{u}\end{bmatrix}=2\re(u^\ast \nabla^2_{z,z}f(z)\, u)+2\re(u^\ast\, 
\nabla^2_{\bar{z},z}f(z)\, \bar{u}).
\end{equation}
In the remaining part of the paper we concentrate on real-valued functions $f$. For minimizing such a function $f$ we shall apply gradient descent 
\begin{equation}\label{eq:Pre1a}
z^t=z^{t-1}-\mu_t\, \nabla f(z^{t-1}) 
\end{equation}
to some appropriate initial vector $z^0\in \C^d$. The parameter $\mu_t > 0$ is called step size. It has to be chosen such that we achieve a descent in every step, viz. $f(z^t)\leq f(z^{t-1})$. The step size can be chosen in different ways. The first option is a constant step size $\mu_t = \mu_c$ for all $t \ge 1$, which is possible to choose, when the action of the Hessian of the of function $f$ is bounded from above as in the next proposition.

\begin{proposition}\label{pro:Pre1}
Let $f:\C^n\to [0,\infty)$ be a twice Wirtinger differentiable function such that the Wirtinger Hessian satisfies 
\begin{equation}\label{eq:Pre3}
\begin{bmatrix} u\\\bar{u}\end{bmatrix}^\ast\, \nabla^2f(z)\, \begin{bmatrix} u\\\bar{u}\end{bmatrix}\leq B\, \left\|\, \begin{bmatrix} u\\\bar{u}\end{bmatrix}\, \right\|_2^2
\end{equation}
for all $z,u\in\C^n$, where $B>0$ is a constant independent of $z$. Let $(z^t)_{t=0}^\infty$ be a sequence generated by \eqref{eq:Pre1a} with arbitrary starting point 
$z^0\in\C^n$ and step size $\mu_t=\mu_c$ such that $0<\mu_c\leq 1/B$ holds.
Then, we have 
\begin{equation}\label{eq:Pre4 f}
f(z^t)- f(z^{t-1})\leq -\mu_t\, \|\nabla_z f(z^t)\|^2_2
\end{equation}
for all $t\geq 1$.\\
In particular, 
\begin{equation}\label{eq:Pre5 f}
\lim_{t\to\infty} \|\nabla_z f(z^t)\|_2^2=0\ \text{ and }\ \min_{t\in\{0,\dots, T\}} \|\nabla_z f(z^t)\|_2^2\leq \frac{f(z^0)}{\mu_c\, (T+1)}.
\end{equation}
\end{proposition}
\banf
The prove \eqref{eq:Pre4 f} is based on the Taylor expansion of $f$ using Wirtinger derivatives, which gives 
$$
f(z+u)=f(z)+\begin{bmatrix} (\partial_z f)^\ast \\ (\partial_{\bar z} f)^\ast \end{bmatrix}^\ast\, \begin{bmatrix} u\\ \bar{u}\end{bmatrix} + 
\begin{bmatrix} u\\ \bar{u}\end{bmatrix}^\ast \int_0^1 (1-s)\nabla^2 f(z+s\, u)\, \dx s\, \begin{bmatrix} u\\ \bar{u}\end{bmatrix},
$$
with $z=z^{t-1}, u=-\mu_c\nabla_z f(z^{t-1})$. Using $\nabla_z f = (\partial_z f)^\ast = \overline{(\partial_{\bar z} f)^\ast}$, inequality \eqref{eq:Pre3}, and the assumption on $\mu_c$ we obtain 
for every $t\geq 1$ 
$$
\begin{array}{ll}
f(z^t)-f(z^{t-1})&\ds\leq -\mu_c\left\|\,  \begin{bmatrix} \nabla_zf(z^{t-1})\\ \overline{\nabla_{z} f(z^{t-1})}\end{bmatrix}\, \right\|^2_2
+\frac{\mu_c^2 B}{2} \left\|\,  \begin{bmatrix} \nabla_zf(z^{t-1})\\ \overline{\nabla_{z} f(z^{t-1})}\end{bmatrix}\, \right\|^2_2\\[2.5ex]
&\leq -2\mu_c(1-\mu_cB/2)\, \|\nabla_z f(z^{t-1})\|_2^2\\[2ex]
&\leq -\mu_c\,  \|\nabla_z f(z^{t-1})\|_2^2. 
\end{array}
$$
To show \eqref{eq:Pre5 f} note that for $T>0$ we have 
$$
\mu_c\, \sum_{t=1}^T \|\nabla_z f(z^{t-1})\|_2^2 \leq f(z^0)-f(z^T)\leq f(z^0), 
$$
which shows in particular that $\sum_{t=1}^\infty \|\nabla_z f(z^{t-1})\|_2^2$ is convergent. Consequently $\|\nabla_z f(z^{t-1})\|_2^2\to 0$ as $t\to\infty$. Finally, we note 
that for $T>0$ 
$$
\min_{t\in\{0,\dots, T\}} \|\nabla_z f(z^t)\|_2^2 \leq \frac{1}{T+1}\, \sum_{t=0}^T \|\nabla_z f(z^{t-1})\|_2^2\leq \frac{f(z^0)}{\mu_c(T+1)}.
$$
\bend

\begin{remark}
We note that Proposition \ref{pro:Pre1} only guarantees convergence to a critical point of the function $f$. This is a common scenario for optimization methods applied to non-convex function $f$. 
\end{remark}

The choice of constant step size in Proposition \ref{pro:Pre1} is based on the worst case scenario for all $z \in \C^{n}$ and may be suboptimal if
$$ 
\begin{bmatrix} u\\\bar{u}\end{bmatrix}^\ast\, \nabla^2f(z)\, \begin{bmatrix} u\\\bar{u}\end{bmatrix} 
\text{ is much smaller than } B \left\|\, \begin{bmatrix} u\\\bar{u}\end{bmatrix}\, \right\|_2^2.
$$
In this case, the so-called Armijo-Goldstein condition can be used to find a bigger step size, for which the decrease of the objective is larger. The Armijo-Goldstein condition 
reads as  
\begin{equation}\label{eq:Pre2}
f(z-\mu\nabla_z f(z))-f(z)\leq - \mu\, \|\nabla_z f(z)\|_2^2.
\end{equation}

A suitable step size is now determined iteratively by the following backtracking line search algorithm, which we will call henceforth Armijo-Goldstein algorithm (AGA, for short). 

\medskip

\begin{algorithm}[H] \label{alg: AGA}
\caption{Backtracking search or Armijo-Goldstein condition (AGA)}
\SetAlgoLined
\SetKwInOut{Input}{Input}
\SetKwInOut{Output}{Output}
\Input{ Differentiable function $f: \C^n \to [0,+\infty)$, current position $z \in \C^n$, initial value $\mu_0>0$, decrease factor $\tau\in (0,1)$.}
\Output{Selected step size $\mu$.}
\For{$j = 0,1,\ldots$}{
\If{$ f(z - \mu_j \nabla_z f(z) ) - f(z) \le - \mu_j \| \nabla_z f(z)\|_2^2 $}
{
\Return{$\mu = \mu_j$\;}
}
$\mu_{j+1} = \tau \mu_j$
}
\end{algorithm}
\medskip
In addition, we make use of the fact that by Proposition \ref{pro:Pre1} the constant step size will always guarantee the desired decrease. Hence, by setting $\mu_0 = \mu_c \tau^{-N}$ for $N \in \N \cup \{0\}$, the AGA will always terminate at $\mu_N = \mu_c$ after $N$ iterations. Note that the so determined parameter $\mu$ depends on $f,z,\tau,\mu_c$ and the number of iterations $N$ to meet the condition, i.e. $\mu=\mu(f,z,\tau,\mu_c,N)$. Moreover, 
\begin{equation}\label{eq: learning rate AGA bounds}
\mu_c \tau^{-N} \geq \mu \geq \mu_c    
\end{equation}
by construction. In case $N=0$, the step size selected by the AGA coincides with the constant step size $\mu_c$. 

\begin{remark}
Often, AGA includes control parameter $0<c<1$, which enters \eqref{eq:Pre2} as an additional multiplier $2c$ on the right hand side. In this paper, we fixed $c=0.5$ to improve readability. However, our results hold true for any $0<c<1$ with a slight adjustment. 
\end{remark}

Regarding the convergence of the gradient descent with $\mu_t$ determined by the AGA we have the following result.

\begin{proposition}\label{pro:Pre1 AGA}
Under conditions of Proposition \ref{pro:Pre1}, a sequence $(z^t)_{t=0}^\infty$ generated by \eqref{eq:Pre1a} with arbitrary starting point 
$z^0\in\C^n$ and step sizes $\mu_t=\mu_t(f,z^{t-1},\tau,\mu_c,N)$ determined by the AGA satisfies \eqref{eq:Pre4 f} and \eqref{eq:Pre5 f}.
\end{proposition}
\banf
If $\mu_t > \mu_c$, inequality \eqref{eq:Pre4 f} holds by construction. Otherwise, for $\mu_t = \mu_c$ Proposition \ref{pro:Pre1} applies. The rest of the proof is analogous to the proof of Proposition \ref{pro:Pre1}.
\bend

\subsection{Reconstruction from quadratic measurements}\label{sec:ReBM}
We will consider in this section the following quadratic reconstruction problem. Suppose we are given data of the form 
\begin{equation}\label{eq:Re0}
y_j=z^\ast Q_j z+\eta_j,\quad j=1,\dots, J,
\end{equation}
where $z\in\C^{dL}$, $Q_j \in \C^{d L \times d L}$ is a positive semidefinite measurement matrix and $\eta_j$ presents the measurement noise. It might look artificial to 
consider $z$ as a vector in  $\C^{d L}$. However, for the polychromatic ptychographic set-up which we will discuss in the next section, $d$ is related to the discretization depth and $L$ is the number of different wavelength. So we have to work with vectors $z
=(z_1,\dots, z_{L})^\top$ consisting of $L$ blocks $z_\ell\in\C^d$.
For the reconstruction of $z$ from data \eqref{eq:Re0} we propose a variant of the so-called Amplitude Flow approach \cite{Xu.2018} and apply a gradient descent for minimizing the related loss function. The loss function we shall consider in this context is 
\begin{equation}\label{eq:Re1}
L_\varepsilon(z)=\sum_{j=1}^J\big[\sqrt{z^\ast Q_j z+\varepsilon}-\sqrt{y_j + \varepsilon}\,\big]^2,\quad z\in\C^{dL},
\end{equation}
where $\varepsilon>0$ is a regularization parameter which is needed to prevent division by zero in the first and second order Wirtinger derivatives. 

We will now add some regularization terms to \eqref{eq:Re1} which are also motivated by the ptychographic imaging application. These are
\begin{itemize}
\item[(a)] Tikhonov regularization 
$$
T(z)= \| z \|_2^2 = \sum_{\ell=1}^{L} \| z_\ell \|^2
$$
\item[(b)] Smoothness. In order to penalize abrupt transition between the different blocks we introduce 
$$
S(z)=\sum_{\ell=1}^{L-1} \kappa_\ell\, \|z_{\ell+1}-z_\ell\|_2^2,
$$
where $\kappa_\ell>0$ are given parameters which we will later associate with the different wavelength in the polychromatic ptychographic set-up. 
\end{itemize}
Putting everything together we arrive at the following regularized loss function 
\begin{equation}\label{eq:Re2}
J(z;\varepsilon, \alpha_T,\alpha_S
)=L_\varepsilon(z)+\alpha_T T(z)+\alpha_S\,S(z)
\end{equation}
with parameters $\varepsilon >0$, $\alpha_T,\alpha_S\ge 0$
. For minimizing \mbox{$J(z):=J(z;\varepsilon, \alpha_T,\alpha_S)
$} for fixed parameters $\epsilon,\alpha_T,\alpha_S
$ we shall apply gradient descent. In order to establish convergence of the gradient descent we make use of Proposition \ref{pro:Pre1 AGA}. For doing so we need   
following auxiliary result. 
\begin{lemma}\label{lem:Re1}
Let $\varepsilon>0$. For the second derivatives of $L_\varepsilon, T, S$, and $R$  the following relations hold 
\begin{align}
\label{eq:Re3a}
\begin{bmatrix} u\\ \bar{u}\end{bmatrix}^\ast\, \nabla^2 L_\varepsilon(z)\, \begin{bmatrix} u\\\bar{u}\end{bmatrix}\leq \Big\|\sum_{j=1}^J\, Q_j\Big\| \
\left\|\,\begin{bmatrix} u\\\bar{u}\end{bmatrix}\,\right\|_2^2,\\
\label{eq:Re3b}
\begin{bmatrix} u\\\bar{u}\end{bmatrix}^\ast\, \nabla^2 T(z)\, \begin{bmatrix} u\\\bar{u}\end{bmatrix}=\left\|\,\begin{bmatrix} u\\\bar{u}\end{bmatrix}\,\right\|_2^2,\\
\label{eq:Re3c}
\begin{bmatrix} u\\\bar{u}\end{bmatrix}^\ast\, \nabla^2 S(z)\, \begin{bmatrix} u\\\bar{u}\end{bmatrix}\leq 
\| K \|
\left\|\,\begin{bmatrix} u\\\bar{u}\end{bmatrix}\,
\right\|_2^2.
\end{align}
for all $z,u\in\C^{dL}$ where the matrix $K\in \R^{L \times L}$ has entries
$$
K_{j,k} = 
\begin{cases}
\kappa_{k-1} (1 - \delta_{k,1} ) + \kappa_{k} (1 - \delta_{k,L}), & j = k, \\
- \kappa_{k}, & j = k+1, \\
- \kappa_{k-1}, & j = k-1, \\
0, & \text{otherwise.}
\end{cases}
$$
In particular, 
\begin{equation}\label{eq:Re3e}
\begin{bmatrix} u\\\bar{u}\end{bmatrix}^\ast\, \nabla^2 J(z)\, \begin{bmatrix} u\\\bar{u}\end{bmatrix}
\leq B \left\|\,\begin{bmatrix} u\\\bar{u}\end{bmatrix}\,\right\|_2^2,
\end{equation}
where $J(z)=J(z;\varepsilon, \alpha_T,\alpha_S)
$ and $B=\Big\|\sum_{j=1}^J\, Q_j\Big\| + \alpha$ with $\alpha:=\alpha_T + \alpha_S \|K\|$.
\end{lemma}

\banf
For the first derivative of $L_\varepsilon$ we obtain 
$$
\nabla_z L_\varepsilon (z)=\sum_{j=1}^J \Big[1-\frac{\sqrt{y_j}}{\sqrt{z^\ast Q_j z +\varepsilon}}\Big]\, Q_j^\ast z.
$$
For the second derivatives of $L_\varepsilon$ we have
$$
\begin{array}{ll}
\nabla_{z,z} L_\varepsilon(z)&=\ds\partial_z\,\sum_{j=1}^J \Big[Q_j^\ast z-\frac{\sqrt{y_j + \varepsilon}\, Q_j^\ast z}{\sqrt{z^\ast Q_j z+\varepsilon}}\Big]\\[2ex]
&=\ds\sum_{j=1}^J\Big[ Q_j^\ast-\frac{\sqrt{y_j + \varepsilon}\, Q_j^\ast}{\sqrt{z^\ast Q_j z+\varepsilon}}+\frac{\sqrt{y_j + \varepsilon}\, Q_j^\ast z\, z^\ast Q_j}{2(z^\ast Q_j z
     +\varepsilon)^{3/2}}\Big],\\[3ex]
\nabla_{\bar{z},z} L_\varepsilon(z)&=\ds\partial_{\bar{z}}\, \sum_{j=1}^J \Big[Q_j^\ast z-\frac{\sqrt{y_j + \varepsilon}\, Q_j^\ast z}{\sqrt{z^\ast Q_j z+\varepsilon}}\Big]
=\ds\sum_{j=1}^J\frac{\sqrt{y_j + \varepsilon}\, Q_j^\ast z\, z^\top Q_j^\top}{2(z^\ast Q_j z+\varepsilon)^{3/2}}.
\end{array}
$$
Moreover, as $L_\varepsilon$ is real-valued, \eqref{eq:Pre2 real val} leads to 
$$
\begin{array}{ll}
\begin{bmatrix} u\\ \bar{u}\end{bmatrix}^\ast\, \nabla^2 L_\varepsilon(z)\, \begin{bmatrix} u\\\bar{u}\end{bmatrix}
=&2\ds\sum_{j=1}^J\Big[u^\ast Q_ju-\frac{\sqrt{y_j + \varepsilon}u^\ast Q_j u}{\sqrt{z^\ast Q_j z+\varepsilon}}+\frac{\sqrt{y_j + \varepsilon}|z^\ast Q_j u|^2}{(z^\ast Q_j z
    +\varepsilon)^{3/2}}\\[2ex]
  &-\ds \frac{\sqrt{y_j + \varepsilon}|z^\ast Q_j u|^2}{2 (z^\ast Q_j z+\varepsilon)^{3/2}}+\frac{\sqrt{y_j + \varepsilon}\, \Re(u^\ast Q_j^* z)^2}{ 2 (z^\ast Q_j z+\varepsilon)^{3/2}}\Big],
\end{array}
$$
where we used that $Q_j$ are Hermitian matrices. Furthermore, note that 
$$
\begin{array}{l}
\ds-\frac{\sqrt{y_j + \varepsilon }u^\ast Q_j u}{\sqrt{z^\ast Q_j z+\varepsilon}}+\frac{\sqrt{y_j + \varepsilon }|z^\ast Q_j u|^2}{(z^\ast Q_j z+\varepsilon)^{3/2}}\\[2ex]
=\ds\frac{\sqrt{y_j + \varepsilon}}{(z^\ast Q_j z+\varepsilon)^{3/2}}\big(|z^\ast Q_j u|^2-u^\ast Q_ju\cdot z^\ast Q_jz\big)-\frac{\varepsilon\sqrt{y_j + \varepsilon}u^\ast Q_j u}
{\sqrt{z^\ast Q_j z+\varepsilon}}\leq 0.
\end{array}
$$
The inequality follows from the fact that $Q_j$ is a positive semidefinite matrix and it can be written as $Q_j=R_j^\ast R_j$, so that
$$
|z^\ast Q_j u|^2\leq\|R_j z\|_2^2\, \|R_j u\|_2^2=z^\ast R_j^\ast R_j z\ u^\ast R_j^\ast R_j u=u^\ast Q_ju\ z^\ast Q_j z.
$$
Moreover, since $\re(u^\ast Q_j z)^2\leq|u^\ast Q_j z|^2$ we arrive at 
$$
\begin{bmatrix} u\\ \bar{u}\end{bmatrix}^\ast\, \nabla^2 L_\varepsilon(z)\, \begin{bmatrix} u\\\bar{u}\end{bmatrix}\leq 2\sum_{j=1}^J u^\ast Q_j u
\leq 2\|Q\|\ \|u\|_2^2= \|Q\|\ \left\|\,\begin{bmatrix} u\\ \bar{u}\end{bmatrix}\,\right\|_2^2,
$$
where $Q=\sum_{j=1}^J Q_j$. \\

For the regularization terms, let us first consider supplementary quadratic function 
$$
H_{M_1,M_2}(z) = z^\ast M_1 z + \re \big( z^\ast M_2 \bar z \big) = z^\ast M_1 z + \tfrac{1}{2}z^\ast M_2 \bar z + \tfrac{1}{2}z^\top \bar M_2 z,
$$
with Hermitian matrices $M_1$ and $M_2$. We note that
$$
\nabla_{z} H_{M_1,M_2}(z) = (z^\ast M_1 + z^\top \bar M_2 )^\ast  = M_1^\ast z + M_2^\top \bar z,
$$
and
$$
\begin{array}{ll}
\nabla_{z,z} H_{M_1,M_2}(z) = M_1, & \nabla_{\bar z,z} H_{M_1,M_2}(z) = M_2^\top = \bar M_2.
\end{array}
$$
Consequently, its Hessian matrix is constant with respect to $z$ and is given by
$$
\nabla^2 H_{M_1,M_2}(z) = 
\begin{bmatrix} 
M_1 & \bar M_2 \\
M_2 & \bar M_1
\end{bmatrix}.
$$
Hence, for $H_{M_1,M_2}$ we have 
\begin{equation}\label{eq: quadratic function bound}
\begin{bmatrix} u\\ \bar{u}\end{bmatrix}^\ast\, \nabla^2 H_{M_1,M_2}(z)\, \begin{bmatrix} u\\\bar{u}\end{bmatrix}
\le \left\| \,
\begin{bmatrix} 
M_1 & \bar M_2 \\
M_2 & \bar M_1
\end{bmatrix} \, \right\|  \left\| \begin{bmatrix} u\\ \bar{u}\end{bmatrix} \right\|_2^2.
\end{equation}
Now, to show \eqref{eq:Re3b} we note that 
$$
T(z) = \| z \|^2 = z^\ast z = H_{I_{dL}, O_{dL}}(z),
$$
where $O_{dL} \in \C^{dL \times dL}$ denotes zero matrix. Equation \eqref{eq: quadratic function bound} yields
$$
\begin{bmatrix} u\\ \bar{u}\end{bmatrix}^\ast\, \nabla^2 H_{I_{dL}, O_{dL}}(z)\, \begin{bmatrix} u\\\bar{u}\end{bmatrix}
\le \left\| \,
\begin{bmatrix} 
I_{dL} & O_{dL} \\
O_{dL} & I_{dL}
\end{bmatrix} \, \right\|  \left\| \begin{bmatrix} u\\ \bar{u}\end{bmatrix} \right\|_2^2
= \left\| \begin{bmatrix} u\\ \bar{u}\end{bmatrix} \right\|_2^2.
$$
To prove \eqref{eq:Re3c} we rewrite
\begin{align*}
S(z) & = \sum_{\ell=1}^{L-1}\kappa_\ell (z_{\ell+1}^\ast-z_\ell^\ast)(z_{\ell+1}-z_\ell) \\
&= \sum_{\ell=1}^{L-1}\kappa_\ell (z_{\ell+1}^\ast z_{\ell+1} + z_\ell^\ast z_\ell - z_\ell^\ast z_{\ell+1} - z_{\ell+1}^\ast z_\ell )
= \sum_{\ell=1}^{L-1} \kappa_\ell z^\ast (K^\ell \otimes I_d) z,
\end{align*} 
where the matrix $K^\ell \in \R^{L \times L}$ has four non-zero entries
$$
K_{\ell,\ell}^\ell = 1, \quad K_{\ell+1,\ell+1}^\ell = 1, \quad K_{\ell,\ell+1}^\ell = -1, \quad K_{\ell+1,\ell}^\ell = -1. 
$$
We recall that for two matrices $M \in \C^{a \times b}, W \in \C^{c \times d}$, the tensor product $M \otimes W \in \C^{ ac \times bd}$ is a block matrix 
\[
M \otimes W 
= 
\begin{bmatrix}
M_{1,1} W & \dots & M_{1,b} W \\
\vdots & \ddots & \vdots \\
M_{a,1} W & \dots & M_{a,b} W \\
\end{bmatrix}.
\]
In the following we will use that the tensor product is linear in both the first and the second arguments and the spectral norm of the tensor product is a product of the spectral norms of its components, that is  $\| M \otimes W\| = \| M \| \cdot \| W \|$.

Returning to $S(z)$, we sum up matrices and apply the above-mentioned linearity of tensor product with respect to the first argument to obtain
$$
S(z) = z^\ast \big( \sum_{\ell=1}^{L-1} \kappa_\ell K^\ell \otimes I_d \big) z = z^\ast (K \otimes I_d) z = H_{K \otimes I_d,O_{dL}}(z).
$$
Hence, by \eqref{eq: quadratic function bound}, we obtain
\begin{align*}
\begin{bmatrix} u\\\bar{u}\end{bmatrix}^\ast\, \nabla^2 S(z)\, \begin{bmatrix} u\\\bar{u}\end{bmatrix}
& \leq 
\left\| \,
\begin{bmatrix} 
K \otimes I_d & O_{dL} \\
O_{dL} & K \otimes I_d
\end{bmatrix} \, \right\|
\left\|\,\begin{bmatrix} u\\\bar{u}\end{bmatrix}\,
\right\|_2^2,
\end{align*}
and by properties of block diagonal matrices and tensor product,
$$
\left\| \,
\begin{bmatrix} 
K \otimes I_d & O_{dL} \\
O_{dL} & K \otimes I_d
\end{bmatrix} \, \right\|
= \| K \otimes I_d \|
= \| K \| \| I_d \|
=\| K \|.
$$

\bend

As an immediate consequence of Proposition \ref{pro:Pre1 AGA} we now obtain the following convergence result. 
\begin{proposition}\label{pro:Re1}
Let $J:\C^{dL}\to [0,\infty)$ be defined as in \eqref{eq:Re2}, so that $J(z)=J(z;\varepsilon, \alpha_T,\alpha_S)
$ with $\varepsilon>0$ and $\alpha_T,\alpha_S \ge 0
$. Consider a sequence $(z^t)_{t=0}^\infty$ generated by \eqref{eq:Pre1a} with an arbitrary starting point $z^0\in\C^{dL}$ and step sizes $\mu_t=\mu_t(J,z^{t-1},\tau,\mu_c,N)$ determined by the AGA. Suppose that the minimal step size satisfies $0<\mu_c\leq 1/B$ where $B$ is the constant in \eqref{eq:Re3e}.
Then we have 
\begin{equation}\label{eq:Pre4}
J(z^t)- J(z^{t-1})\leq-\mu_t\, \|\nabla_z J(z^t)\|^2_2
\end{equation}
for all $t\geq 1$.\\
In particular, 
\begin{equation}\label{eq:Pre5}
\lim_{t\to\infty} \|\nabla_z J(z^t)\|_2^2=0\ \text{ and }\ \min_{t\in\{0,\dots, T\}} \|\nabla_z J(z^t)\|_2^2\leq \frac{J(z^0)}{\mu_c\, (T+1)}.
\end{equation}
\end{proposition}

\begin{Remark}
We note that condition $\varepsilon >0$ can be relaxed to $\varepsilon \ge 0$. For the case $\varepsilon =0$, the loss function $L_0(z)$ is not everywhere differentiable. This can be solved by considering the generalized gradient $\nabla L_0(z) = \lim_{\varepsilon \to 0+} \nabla L_\varepsilon(z)$ instead of standard gradient. In this case the analogue of Proposition \ref{pro:Re1} holds. For details, see \cite{Xu.2018}. 
\end{Remark}

\section{Polychromatic Ptychography}\label{sec:PP}
In this section we will address the problem of reconstructing an object from  measurements \eqref{eq: measurements continuous}. More precisely, the measurements are first discretized and, then, recovery of discretized object and the window is considered. 


\subsection{Discretization of problem}

For discretization of the single Fourier integrals in \eqref{eq: measurements continuous} we are using the grid $\Gamma_d=\frac{1}{d}\,  [d]^2=\frac{1}{d}\, [d]\times [d]$, where we have used the notation $[d]=\{0,1,\dots, d-1\}$, and we 
consider shifts $r\in\Gamma_d$. Also, without loss of generality we may assume that all $f_{\lambda_\ell}$ are supported in $\mathcal{Q}=[0,1]^2$. This gives 
\begin{equation*}
\mathscr{F}(f_{\lambda_\ell}\, T_{m/d}w_{\lambda_\ell})(\xi/\lambda_\ell)\approx\frac{1}{d^2}\sum_{n\in[d]^2}f_{\lambda_\ell}({\ts\frac{n}{d}})\, 
{w}_{\lambda_\ell}({\ts\frac{n-m}{d}}) \mathbbm{1}_{\mathcal{Q}}(\ts\frac{n-m}{d})\, \e^{-2\pi\im \frac{1}{\lambda_\ell}\xi\cdot \frac{1}{d}n},
\end{equation*}
for $m \in \mathcal M$, with $\mathcal M \subseteq [-d,d]^2$ denoting a set of observed shift position.
Note that we are not working with cyclic shifts of the mask but cutting out that 
part of $w_{\lambda_\ell}$ which lies in $\mathcal{Q}$. \\ 
The dual grid for the Fourier transform dilated by $1/\lambda_\ell$ is $\hat{\Gamma}_{\lambda_\ell}=\lambda_\ell\,[d]^2$ and in order to avoid multiple 
contributions from the smallest wavelength term in \eqref{eq: measurements continuous} we have to evaluate intensity function on the dual grid $\hat{\Gamma}_{\lambda_1}$. 
\begin{figure}
    \centering
    \includegraphics
    [width = 1\linewidth ]
    {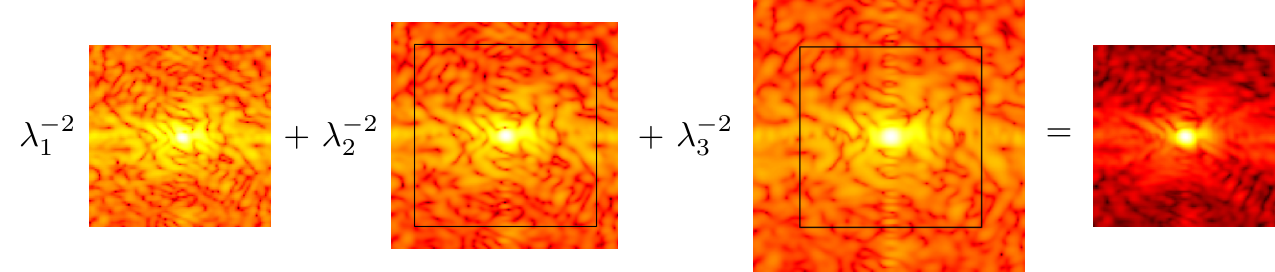}
    \caption{Example of diffraction patterns resulted from polychromatic illumination with $d = 100\times100$, $L=3$, $\lambda = (1,1.25,1.5)$.
    The black rectangles denote the dual grid $\hat{\Gamma}_{\lambda_1}$.}
    \label{fig:2}
\end{figure}

Using the 
notation $x_\ell=\big(f_{\lambda_\ell}({\ts\frac{n}{d}})\big)_n$ and $w_\ell = \big( w_{\lambda_\ell}({\ts\frac{n}{d}}) \big)_n$, 
we get the discretized multi-spectral intensity measurements 
\begin{align}
\mathfrak{I}(m/d,k\lambda_1) & \approx y_{m,k} + \eta_{m,k} \label{eq: poly meas} \\
& := \sum_{\ell=1}^L \frac{1}{\lambda_\ell^2} ~
\big| \sum_{n\in[d]^2} (x_\ell)_{n} (w_\ell)_{n-m} \mathbbm{1}_{\mathcal{Q}}(\ts\frac{n-m}{d}) \e^{-2\pi\im \frac{\lambda_1}
{\lambda_\ell}\,k\cdot \frac{1}{d}n} \, \big|^2 + \eta_{m,k}, \nonumber
\end{align}
where $\eta$ denotes composite noise from the discretization process and from the measurement process.
\subsection{Non-blind problem}\label{sec:NBP}

First, let us consider the recovery of the unknown vector $x$ under the assumption that $w$ is known. Then, we can rewrite the measurements as
\begin{equation}\label{eq:PP4}
y_{m,k}= \sum_{\ell=1}^L \frac{1}{\lambda_\ell^2} \, \left|\, \langle x_{\ell},  a_{m,k}^{\ell}\rangle\, \right|^2,\qquad m,k\in [d]^2,
\end{equation}
with the vectors $a_{m,k}^{\ell}=\big( (\overline{w}_{\ell})_{n-m} \mathbbm{1}_{\mathcal{Q}}(\ts\frac{n-m}{d}) \e^{2\pi\im \frac{\lambda_1}
{\lambda_\ell}\,k\cdot \frac{1}{d}n}\big)_n$.
Accordingly, we have $y=(y_{m,k})_{m,k}$. Relation \eqref{eq:PP4} can also be written in a bilinear form 
\begin{align}
y_{m,k} & = \sum_{\ell=1}^L \frac{1}{\lambda_\ell^2} \, \overline{\langle x_{\ell},  a_{m,k}^{\ell}\rangle}   \langle x_{\ell},  a_{m,k}^{\ell}\rangle
= \sum_{\ell=1}^L \frac{1}{\lambda_\ell^2} \, \langle a_{m,k}^{\ell}, x_{\ell}\rangle\langle x_{\ell},  a_{m,k}^{\ell}\rangle 
\nonumber \\
& = \sum_{\ell=1}^L \frac{1}{\lambda_\ell^2} x_{\ell}^\ast a_{m,k}^{\ell} (a_{m,k}^{\ell})^* x_{\ell}
=z^\ast\, Q_{m,k}^w z \label{eq:PP meas}
\end{align}
with $z=(x_1,\dots, x_L)^\top\in\C^{dL}$ and the block diagonal matrix $Q_{m,k}^w \in\C^{dL\times dL}$ with the diagonal blocks 
$\lambda_\ell^{-2} \, a_{m,k}^\ell(a_{m,k}^\ell)^\ast$. We use an upper index $w$ to emphasize on dependence of $Q_{m,k}^w$ on the mask $w$. 
In the form \eqref{eq:PP5} the polychromatic ptychographic reconstruction problem is now a problem of the form \eqref{eq:Re0} which was considered in the previous section and it can be solved by minimizing $J$ as given in \eqref{eq:Re2}. The smoothness penalty $S(z)$ imposes the continuity with respect to $\lambda$ of the object. 
The purpose of the Tikhonov regularization $T(z)$ will be explained later, when we turn to the reconstruction of both the object and the window.

For our discussion we shall make use of a matrix notation of the expression \eqref{eq:PP4}, viz. 
\begin{equation}\label{eq:PP5}
y=\sum_{\ell=1}^L \frac{1}{\lambda_\ell^2} \, |A_\ell\, x_\ell |^2,
\end{equation}
where  $A_\ell=(A_{m}^{\ell})_{m}$ is a row block matrix with row blocks $A_{m}^{\ell}$ with rows $(a_{m,k}^{\ell})_{k}$. Note that every row block corresponds to one specific shift of the mask and is given by
$$
A_{m}^{\ell} = F^\ell D^\ell_m,
$$
with the matrix $F^\ell$ defined as
\begin{equation} \label{eq: F ell matrix}
(F^\ell)_{k,n} = \e^{-2\pi\im \frac{\lambda_1} {\lambda_\ell}\,k\cdot \frac{1}{d}n},
\end{equation}
and diagonal matrix $D^\ell_m$ given by
$$
(D^\ell_m)_{n,n} = (w_{\ell})_{n-m} \mathbbm{1}_{\mathcal{Q}}(\ts\frac{n-m}{d}).
$$

Consequently, we obtain the following corollary to Proposition \ref{pro:Re1}.
\begin{Corollary}\label{col: non-blind ptycho}
Consider measurements of the form \eqref{eq: poly meas}. Let $J:\C^{dL}\to [0,\infty)$, $J(z)=J(z;\varepsilon, \alpha_T,\alpha_S
)$ as defined in \eqref{eq:Re2} with matrices $Q_{m,k}^w$ as in \eqref{eq:PP meas}. Suppose that for the minimal step size satisfies $0<\mu_c\leq 1/B(w)$ with $B(w)$ given by
\begin{equation}\label{eq: bound for object}
B(w) = \max_{\ell = 1, \dots, L} \left\{ \frac{\| F_\ell \|^2}{\lambda_\ell^2}  \max_{n \in [d]} \big[ \sum_{m \in \mathcal M} |(w_\ell)_{n-m}|^2 \mathbbm{1}_{\mathcal{Q}}(\ts\frac{n-m}{d})  \big] \right\} + \alpha,
\end{equation}
where $F_\ell$ are matrices defined in $\eqref{eq: F ell matrix}$ and $\alpha$ as in Lemma \ref{lem:Re1}. 
Then, the results of Proposition \ref{pro:Re1} apply to the sequence $(z^t)_{t=0}^\infty$ generated by \eqref{eq:Pre1a} with an arbitrary starting point $z^0\in\C^{dL}$ and step sizes $\mu_t=\mu_t(J,z^{t-1},\tau,\mu_c,N)$ determined by the AGA.
\end{Corollary}
\banf
We only need show that $\mu_c$ satisfies condition in Proposition \ref{pro:Re1}, which is equivalent to proving that $B(w)$ in the statement of the theorem is greater or equal than $B$ in \eqref{eq:Re3e}. Furthermore, since all but one summands are the same, we only need to show the inequality
$$
\big\| \sum_{m \in \mathcal M}\sum_{k \in [d]} Q_{m,k}^w \big\| \leq \max_{\ell = 1, \dots, L} \left\{ \frac{\| F_\ell \|^2 }{\lambda_\ell^2} \max_{n \in [d]} \big[ \sum_{m \in \mathcal M} |(w_\ell)_{n-m}|^2 \mathbbm{1}_{\mathcal{Q}}(\ts\frac{n-m}{d})  \big] \right\}.
$$
Let us first compute $\big\| \ds\sum_{m \in \mathcal M} \ds\sum_{k \in [d]} Q_{m,k}^2 \big\|$. Since each $Q_{m,k}^w$ is block diagonal,  the sum of $Q_{m,k}^w$ is block diagonal as well. Moreover, the spectral norm of the block diagonal matrix is the maximum of the spectral norms of all blocks. Thus, we have
\begin{align}
\big\| \sum_{m \in \mathcal M}\sum_{k \in [d]} Q_{m,k}^w \big\|
& = \max_{\ell = 1, \dots, L} \big\| \sum_{m \in \mathcal M}\sum_{k \in [d]} \frac{1}{\lambda_\ell^2} \, a_{m,k}^\ell(a_{m,k}^\ell)^\ast \big\| \nonumber \\
& = \max_{\ell = 1, \dots, L} \frac{1}{\lambda_\ell^2} \big\| \sum_{m \in \mathcal M} (A_{m}^{\ell})^\ast A_{m}^{\ell} \big\| \nonumber \\
& = \max_{\ell = 1, \dots, L} \lambda_\ell^{-2} \| (A_\ell)^\ast A_\ell \|  = \max_{\ell = 1, \dots, L}\lambda_\ell^{-2} \| A_\ell \|^2. \label{eq: spec norm tech 1}
\end{align}
Matrix $A_\ell$ is a row block matrix with blocks $A_{m}^{\ell} = F^\ell D^\ell_m$. Therefore, we can presents $A_\ell = \tilde F_\ell D_\ell$ as a product of block diagonal matrix $\tilde F_\ell$ with blocks $F^\ell$ and row block matrix $D_\ell =(D_{m}^{\ell})_{m}$.
Hence, 
\begin{equation}\label{eq: spec norm tech 2}
\| A_\ell \| = \| \tilde F_\ell D_\ell \| \leq \| \tilde F_\ell \| \| D_\ell \|
\end{equation}
The matrix $\tilde F_\ell$ is block diagonal with blocks $F^\ell$ and thus $\| \tilde F_\ell \| = \| F^\ell \|$. For $\| D_\ell \|$, we observe that $\| D_\ell \|^2 = \| D_\ell^\ast D_\ell \|$ and 
$$
D_\ell^\ast D_\ell
= \sum_{m \in \mathcal M} (D^\ell_m)^\ast D^\ell_m .
$$
Recalling that $D^\ell_m$ are diagonal yields that $D_\ell^\ast D_\ell$ is again diagonal and its entries are given by
\begin{align}
(D_\ell^\ast D_\ell)_{n,n} & = \big( \sum_{m \in \mathcal M} (D^\ell_m)^\ast D^\ell_m \big)_{n,n} \nonumber 
= \sum_{m \in \mathcal M} \overline{(D^\ell_m)}_{n,n} (D^\ell_m)_{n,n} \nonumber \\
& = \sum_{m \in \mathcal M} |(w_\ell)_{n-m}|^2 \mathbbm{1}_{\mathcal{Q}}(\ts\frac{n-m}{d}). \label{eq: spec norm tech 3}
\end{align}
Finally, we combine \eqref{eq: spec norm tech 1}, \eqref{eq: spec norm tech 2} and \eqref{eq: spec norm tech 3} to obtained desired result,
\begin{align*}
\big\| \sum_{m \in \mathcal M}\sum_{k \in [d]} Q_{m,k}^w \big\| 
& \leq \max_{\ell = 1, \dots, L} \lambda_\ell^{-2} \| \tilde F_\ell \|^2 \| D_\ell \|^2
= \max_{\ell = 1, \dots, L} \lambda_\ell^{-2} \| F_\ell \|^2 \| D_\ell^\ast D_\ell \| \\
& = \max_{\ell = 1, \dots, L} \lambda_\ell^{-2} \| F_\ell \|^2 \max_{n \in [d]} \big[ \sum_{m \in \mathcal M} |(w_\ell)_{n-m}|^2  \mathbbm{1}_{\mathcal{Q}}(\ts\frac{n-m}{d}) \big].
\end{align*}
\bend

Now let us consider a dual recovery problem of finding the mask under assumption that the object is known.
Changing the summation order from $n$ to $n+m$ in \eqref{eq: poly meas} and keeping in mind that $f_{\lambda_\ell}$ is only supported on $\mathcal{Q}$ we arrive at
\begin{align*}
y_{m,k} & = \sum_{\ell=1}^L \frac{1}{\lambda_\ell^2} \,
\big| \sum_{n\in[d]^2} (x_\ell)_{n + m} (w_\ell)_{n} \mathbbm{1}_{\mathcal{Q}}(\ts\frac{n+m}{d}) \e^{-2\pi\im \frac{\lambda_1}
{\lambda_\ell}\,k\cdot \frac{1}{d} (n + m)} \big|^2\\
& = \sum_{\ell=1}^L \frac{1}{\lambda_\ell^2} \,
\big| \sum_{n\in[d]^2} (x_\ell)_{n + m} (w_\ell)_{n} \mathbbm{1}_{\mathcal{Q}}(\ts\frac{n+m}{d}) \e^{-2\pi\im \frac{\lambda_1}
{\lambda_\ell}\,k\cdot \frac{1}{d} n} \big|^2.
\end{align*}
We can rewrite the measurements as
\begin{equation}\label{eq: meas window}
y_{m,k}=w^\ast\, Q_{m,k}^x w
\end{equation}
with $w=(w_1,\dots, w_L)^\top\in\C^{dL}$ and the block diagonal matrix $Q_{m,k}^x \in\C^{dL\times dL}$ with rank-one diagonal blocks 
$\lambda_\ell^{-2} \, v_{m,k}^\ell(v_{m,k}^\ell)^\ast \in \C^{d \times d}$ formed by vectors $v_{m,k}^{\ell} = \big( (\bar x_{\ell})_{n+m} \mathbbm{1}_{\mathcal{Q}}(\ts\frac{n+m}{d}) \e^{2\pi\im \frac{\lambda_1}
{\lambda_\ell}\,k\cdot \frac{1}{d}n}\big)_n$.

Analogously to Corollary \ref{col: non-blind ptycho} we obtain the following.
\begin{Corollary}\label{col: non-blind ptycho window}
Consider measurements of the form \eqref{eq: poly meas}. Let $J:\C^{dL}\to [0,\infty)$, $J(w)=J(w;\varepsilon, \beta_T,\beta_S
)$ as defined in \eqref{eq:Re2} with matrices $Q_{m,k}^x$ as in \eqref{eq: meas window}. Suppose that for the minimal step size satisfies $0<\nu_c\leq 1/B(x)$ with $B(x)$ given by
\begin{equation}\label{eq: bound for window}
B(x) = \max_{\ell = 1, \dots, L} \left\{ \frac{ \| F_\ell \|^2}{\lambda_\ell^2} \max_{n \in [d]} \big[ \sum_{m \in \mathcal M} |(x_\ell)_{n+m}|^2 \mathbbm{1}_{\mathcal{Q}}(\ts\frac{n+m}{d}) \big] \right\} + \beta,
\end{equation}
where $F_\ell$ are matrices defined in $\eqref{eq: F ell matrix}$ and $\beta = \beta_T + \beta_S \| K \|$ defined analogously to $\alpha$ in Lemma \ref{lem:Re1}. 
Then, the results of Proposition \ref{pro:Re1} apply to the sequence $(w^t)_{t=0}^\infty$ generated by \eqref{eq:Pre1a} with an arbitrary starting point $w^0\in\C^{dL}$ and step sizes $\nu_t=\nu_t(J,w^{t-1},\tau,\nu_c,N)$ determined by the AGA.
\end{Corollary}

\subsection{Blind problem}\label{sec:BP}

Finally, we turn to the reconstruction of both $z$ and $w$ from the measurements \eqref{eq: poly meas}.
The problem can now be formulated as 
$$
y_{m,k} + \eta_{m,k}= z^\ast Q_{m,k}^{w} z +\eta_{m,k} = w^\ast Q_{m,k}^{z} w +\eta_{m,k},
$$
with $Q_{m,k}^{w}$ as in \eqref{eq:PP meas} and $Q_{m,k}^{z}$ as in \eqref{eq: meas window}. 
The modified regularized loss function for blind polychromatic ptychography takes the form 
\begin{equation}\label{eq:PP8}
\J(z,w;\varepsilon, \alpha_T,\alpha_S,
\beta_T, \beta_S
)
=L_\varepsilon(z;w)+\alpha_T T(z)+\alpha_S\,S(z)
+\beta_T T(w)+\beta_S S(w)
\end{equation}
with parameters $\varepsilon>0, \alpha_T, \alpha_S,
\beta_T,\beta_S
\ge 0$ and $L_\varepsilon(z;w)$ denoting $L_\varepsilon(z)$ with matrices $Q_{m,k}^{w}$. 
The reconstruction process can now be described by the following iterative method. Note that 
\begin{align}
& \J(z,w)\big|_{w = const.} = J(z; \alpha_T,\alpha_S,\alpha_R) + C_w, \label{eq: BJ to J relation}\\
& \J(z,w)\big|_{z = const.} = J(w; \beta_T,\beta_S,\beta_R) + C_z \nonumber ,
\end{align}
where $J$ is defined as in \eqref{eq:Re2} and constants $C_z$ and $C_w$ which are depending on $z$ and $w$, respectively. Hence, the partial gradients $\nabla_z \J$ resp. $\nabla_w \J$ are given by 
\begin{align}
\nabla_z \J(z,w) & = \nabla_z J(z;\alpha_T,\alpha_S,\alpha_R) \label{eq:PP9} \\
\nabla_w \J(z,w) & = \nabla_w J(w; \beta_T,\beta_S,\beta_R) \nonumber
\end{align}

Naturally, the reasonable choice is to perform the full gradient descent. With initial guesses $z^0, w^0$, the gradient descent iterations 
$$
\begin{bmatrix} z^{t}\\ w^{t} \end{bmatrix} = \begin{bmatrix} z^{t-1}\\ w^{t-1} \end{bmatrix} 
- \mu_t 
\begin{bmatrix} \nabla_z \J(z^{t-1},w^{t-1}) \\
\nabla_w \J(z^{t-1},w^{t-1})
\end{bmatrix}
$$
can be performed. However, the difficulty behind this approach is to find a suitable step size $\mu_t$, which depends on both $z^{t-1}$ and $w^{t-1}$. Instead, we consider an alternative approach reminiscent of the alternating minimization \cite{Beck.2015}. That is, with initial guesses $z^0, w^0$, we fix the current mask iterate $w^t$ and perform a fixed number $I_z \in \N$ gradient descent iterations for the object, which provides the new iterate $z^{t+1}$. Next, we fix $z^{t+1}$ and perform $I_w \in \N$ iterations with respect to the mask resulting in $w^{t+1}$.
The sequences $\mu_{t,i}$ and $\nu_{t,j}$ are appropriate step sizes determined by the AGA as discussed in Section \ref{sec:Pre}.
Note that due to fixation of the variables, the constant step sizes can be chosen accordingly to Corollaries \ref{col: non-blind ptycho} and \ref{col: non-blind ptycho window}.

This results in the following procedure.


\medskip

\begin{algorithm}[H] \label{alg:PP1}
\caption{Amplitude flow for blind polychromatic ptychography}
\SetAlgoLined
\SetKwInOut{Input}{Input}
\SetKwInOut{Output}{Output}
\Input{Measurements $y$ as in \eqref{eq: poly meas}, starting points $z^0,w^0 \in \C^{dL}$, number of iterations $T \in \N$, number of object and mask iterations $I_z, I_w \in \N$, regularization parameters $\varepsilon>0$, $\alpha_T, \alpha_S
\beta_T, \beta_S
\ge 0$, AGA parameters $N \in \N, 0<\tau<1.$}
\Output{$z,w\in \C^{dL}$}
\For{$t = 1,\ldots, T$}{
Let $z^{t,0} = z^{t-1}$\;
Set $\mu_{t,c} = 1/B(w^{t-1})$\;
\For{$i = 1, \ldots, I_z$}
{
Select $\mu_{t,i} = \mu_{t,i}(\, \J\,\big|_{w = w^{t-1}}, z^{t,i-1},\tau,\mu_{t,c},N)$ via AGA\;
Update $z^{t,i} = z^{t,i-1} - \mu_{t,i} \nabla_z \J( z^{t,i-1}, w^{t-1})$\;
}
Let $z^t = z^{t,I_z}$ and $w^{t,0} = w^{t-1}$ \;
Set $\nu_{t,c} = 1/B(z^{t})$\;
\For{$j = 1, \ldots, I_w$}
{
Select $\nu_{t,j} = \nu_{t,j}(\, \J\, \big|_{z = z^{t}},w^{t,j-1},\tau,\nu_{t,c},N)$ via AGA\;
Update $w^{t,j} = w^{t,j-1} - \nu_{t,j} \nabla_w \J( z^{t}, w^{t,j-1})$\;
}
Let $w^t = w^{t,I_w}$\;
}
\Return{$z = z^{T}, w = w^{T}$}
\end{algorithm}

\medskip

First, let us show that the step sizes $\mu_{t,i}$ and $\nu_{t,j}$ are always finite.
\begin{lemma}\label{l: learning rates finite}
Assume that $\varepsilon,\alpha_T, \beta_T>0$. Then, for all $z,w \in \C^{dL}$ we have
$B(w) \ge \alpha_T > 0$ and $B(z) \ge \beta_T > 0$.
Furthermore, step sizes $\mu_{t,i}$ and $\nu_{t,j}$, $i = 1, \dots, I_z$, $j = 1, \dots, I_w$, $t \ge 1$, determined by Algorithm \ref{alg:PP1} are bounded by $\tau^{-N}/\alpha_T$ and $\tau^{-N}/\beta_T$, respectively.
\end{lemma}
\banf
By definition \eqref{eq: bound for object}, $B(w)$ satisfies $B(w) \ge \alpha_T >0$ and, therefore, for all step sizes $\mu_{t,i}$ selected via AGA, by \eqref{eq: learning rate AGA bounds}, we have
$$
\mu_{t,i} \le \tau^{-N} \mu_{t,c} = \tau^{-N}/B(w^{t-1,I_w}) \le \tau^{-N}/ \alpha_T < \infty.
$$
Analogously, $\nu_{t,j} \le \tau^{-N}/\beta_T < \infty$. 
\bend


Note that bounds $B(z^t)$ and $B(w^t)$ are not constant and, hence, theory established in \cite{Beck.2015} is not applicable. Thus, we derive the following results regarding the convergence of proposed alternating minimization process.
\begin{theorem}\label{theo:PP1}
Let $\J:\C^{dL}\times\C^{dL}\to [0,\infty)$ be defined as in \eqref{eq:PP8} with $\varepsilon, \alpha_T, \beta_T >0$ and assume the two sequences $(z^t)_{t\geq 0}, (w^t)_{t\geq 0}$ are generated by 
Algorithm \ref{alg:PP1} with arbitrary starting points $z^0, w^0\in\C^{dL}$ and let $\mu_{t,i}$ and $\nu_{t,j}$ be step sizes as determined by Algorithm \ref{alg:PP1}.
Then for each subiteration of Algorithm \ref{alg:PP1} we have 
\begin{align}
& \J(z^{t,i},w^{t-1})-\J(z^{t,i-1},w^{t-1})\leq -\mu_{t,i}\, \|\nabla_z \J(z^{t,i-1},w^{t-1})\|_2^2 \label{eq: object subit}\\
& \J(z^{t},w^{t,j})-\J(z^{t},w^{t,j-1})\leq -\nu_{t,j}\, \|\nabla_w \J(z^{t},w^{t,j-1})\|_2^2. \label{eq: window subit} 
\end{align}
for every $t\geq1$ and $i=1,\dots, I_z,\ j=1,\dots, I_w$. \\

Moreover, 
$$
\lim_{t\to\infty} \|\nabla \J(z^{t},w^{t})\|_2^2=0,
$$
where the rate of convergence is dominated by the number 
$$
\frac{\max\{ \alpha_T^{-1}, \beta_T^{-1} \} (\J(z^0,w^0))^2 \max_{\ell = 1, \dots, L} \lambda_\ell^{-2} \| F_\ell \|^2 + 
\J(z^0,w^0) \max\{ \alpha, \beta\}
}{T \min\{I_z, I_w \} },
$$
with $\alpha$ and $\beta$ as in Corollaries \ref{col: non-blind ptycho} and \ref{col: non-blind ptycho window}.
\end{theorem}
\banf
Let $t\geq1$ be fixed. For each the object subiteration, $w^{t-1}$ is fixed and the constant step size $\mu_{t,c} = 1/B(w^{t-1})$ satisfies the conditions of Corollary \ref{col: non-blind ptycho}. Thus, we have 
\begin{align*}
J(z^{t,i}; \varepsilon,\alpha_T,\alpha_S,\alpha_R)-J(z^{t,i-1}; \varepsilon,\alpha_T,\alpha_S,\alpha_R) \\
\quad \quad \quad \quad \leq -\mu_{t,i}\, \|\nabla J(z^{t,i-1}; \varepsilon,\alpha_T,\alpha_S,\alpha_R)\|_2^2.
\end{align*}
In a view of equalities \eqref{eq: BJ to J relation} and \eqref{eq:PP9}, we obtain \eqref{eq: object subit}. Analogously, Corollary \ref{col: non-blind ptycho window} yields \eqref{eq: window subit}. These estimates show that every subiteration of the Algorithm \ref{alg:PP1} reduces the value of the loss function $\J(z,w)$.
Furthermore,
$$
\J(z^{t,I_z},w^{t-1})-\J(z^{t,0},w^{t-1})\leq -\sum_{i=1}^{I_z}\mu_{t,i} \|\nabla_z \J(z^{t,i-1},w^{t-1})\|_2^2
$$
and analogously we get 
$$
\J(z^{t},w^{t,I_w})-\J(z^{t},w^{t,0})\leq -\sum_{j=1}^{I_w}\nu_{t,j} \|\nabla_w \J(z^{t},w^{t,j-1})\|_2^2.
$$
Note that by construction $z^{t,I_z} = z^t$, $z^{t,0} = z^{t-1}$, $w^{t,I_w} = w^t$, $w^{t,0} = w^{t-1}$.
Hence, combining the inequalities leads to
$$
\sum_{i=1}^{I_z}\mu_{t,i} \|\nabla_z \J(z^{t,i-1},w^{t-1})\|_2^2+\sum_{j=1}^{I_w}\nu_{t,j} \|\nabla_w \J(z^{t},w^{t,j-1})\|_2^2
\leq \J(z^{t-1},w^{t-1}) - \J(z^{t},w^{t}). 
$$
For a fixed $T\geq 1$ a summation over all $t=1,\dots, T$ gives 
\begin{align}
&\sum_{t=1}^T\Big[\sum_{i=1}^{I_z}\mu_{t,i} \|\nabla_z \J(z^{t,i-1},w^{t-1})\|_2^2+\sum_{j=1}^{I_w}\nu_{t,j} \|\nabla_w \J(z^{t},w^{t,j-1})\|_2^2\Big] \nonumber \\
& \quad \quad \leq \sum_{t =1}^{T} \left[ \J(z^{t-1},w^{t-1})-\J(z^{t},w^t) \right] \nonumber \\
& \quad \quad = \J(z^{0},w^{0}) - \J(z^{T},w^T) \label{eq: objective decreas} \\
& \quad \quad \leq  \J(z^0,w^0), \nonumber
\end{align}
where we used that $\J(z,w) \ge 0$. Hence, by taking $T \to \infty$, we arrive at 
$$
\sum_{t=1}^\infty\Big[\sum_{i=1}^{I_z}\mu_{t,i} \|\nabla_z \J(z^{t,i-1},w^{t-1})\|_2^2+\sum_{j=1}^{I_w}\nu_{t,j} \|\nabla_w \J(z^{t},w^{t,j-1})\|_2^2\Big]<\infty,
$$
which implies that 
$$
\sum_{i=1}^{I_z}\mu_{t,i} \|\nabla_z \J(z^{t,i-1},w^{t-1})\|_2^2+\sum_{j=1}^{I_w}\nu_{t,j} \|\nabla_w \J(z^{t},w^{t,j-1})\|_2^2\to 0
$$
as $t\to\infty$. Since all terms are non-negative we eventually get  
$$
\mu_{t,i} \|\nabla_z \J(z^{t,i-1},w^{t-1})\|_2^2\to 0,\quad\text{and}\quad  \nu_{t,j} \|\nabla_w \J(z^{t},w^{t,j-1})\|_2^2\to 0
$$ 
for all $i=1,\dots, I_z,\ j=1,\dots, I_w$. 

In order to show the desired convergence for the norm of the gradients we have to show that the step sizes $\mu_{t,i}$ and $\nu_{t,j}$ are not converging to zero as $t\to\infty$ for all $i,j$. Since the step size $\mu_{t,i}$ is assumed to be determined by the AGA we have 
\begin{equation}\label{eq: learinign rate below bound}
\mu_{t,i} \geq
\mu_{t,c} = 1/B(w^{t-1}),
\end{equation}
where $\mu_{t,c}$ is the minimal step size. Hence, showing that $\mu_{t,i}$ does not vanish is equivalent to prove that the sequence $(B(w^{t}))_{t \ge 0}$ is bounded from above. Recall that by Corollary \ref{col: non-blind ptycho}, $B(w^{t})$ is given by
\begin{align*}
B(w^{t}) & = \max_{\ell = 1, \dots, L} \left\{ \lambda_\ell^{-2} \| F_\ell \|^2 \max_{n \in [d]} \big[ \sum_{m \in \mathcal M} |(w_\ell^{t})_{n-m}|^2 \mathbbm{1}_{\mathcal{Q}}(\ts\frac{n-m}{d})  \big] \right\} + \alpha\\
& \le \max_{\ell = 1, \dots, L} \left\{ \lambda_\ell^{-2} \| F_\ell \|^2 \|w_\ell^{t} \|_2^2 \right\} + \alpha.
\end{align*}
Consequently, $( B(w^{t}) )_{t \ge 0}$ is bounded if and only if all sequences $( \|w_\ell^{t} \|_2^2 )_{t \ge 0}$, $\ell = 1, \dots, L$ are bounded. Let us show by contradiction that $( \|w_\ell^{t} \|_2^2 )_{t \ge 0}$, $\ell = 1, \dots, L$ are bounded by $\J(z^0,w^0) / \beta_T$. More precisely, assume that for some $\ell_0 \in \{1,\dots, L\}$, the sequence $( \|w_{\ell_0}^{t} \|_2^2 )_{t \ge 0}$ exceeds $\J(z^0,w^0) / \beta_T$. Then, there exists $t_0 \ge 0$ such that $\|w_{\ell_0}^{t_0} \|_2^2 > \J(z^0,w^0) / \beta_T$ and we obtain
$$
\J(z^{t_0}, w^{t_0}) \ge \beta_T \sum_{\ell = 1}^L \| w_{\ell}^{t_0} \|_2^2 \ge \beta_T \| w_{\ell_0}^{t_0} \|_2^2 > \J(z^0,w^0),
$$
which is not possible, since we showed in \eqref{eq: objective decreas} that with each iteration objective does not increase. Therefore, all $( \|w_\ell^{t} \|_2^2 )_{t \ge 0}$ are bounded from above by $\J(z^0,w^0) / \beta_T$ and $( B(w^{t}) )_{t \ge 0}$ is also bounded from above by
\begin{equation} \label{eq: upper bound window}
B_w :=  \beta_T^{-1} \J(z^0,w^0) \max_{\ell = 1, \dots, L} \lambda_\ell^{-2} \| F_\ell \|^2  + \alpha >0, 
\end{equation}
where the strict inequality is valid since $B_w \ge \alpha_T > 0$. Hence, by \eqref{eq: learinign rate below bound}, 
\begin{equation} \label{eq: lower bound learning rate}
\mu_{t,i} \ge 1/B_w > 0,
\end{equation}
and we obtain,
$$
\|\nabla_z \J(z^{t,i-1},w^{t-1})\|_2^2 \to 0,
$$
as $t \to \infty$ for all $i = 1, \dots, I_z$. Similarly, $( B(z^{t}) )_{t \ge 0}$ is bounded from above by 
$$
B_z := \alpha_T^{-1} \J(z^0,w^0) \max_{\ell = 1, \dots, L} \lambda_\ell^{-2} \| F_\ell \|^2 + \beta >0,
$$ 
so that for all $j =1, \dots, I_w$, we obtain 
$$
\nu_{t,j} \ge 1/B_z > 0
$$
and
$$
\|\nabla_w \J(z^{t},w^{t,j})\|_2^2 \to 0, 
$$
as $t \to \infty$. In particular,
\begin{align*}
& \|\nabla_z \J(z^{t},w^{t})\|_2^2 = \|\nabla_z \J(z^{t+1,0},w^{t})\|_2^2 \to 0,\\
& \|\nabla_w \J(z^{t},w^{t-1})\|_2^2 = \|\nabla_w \J(z^{t},w^{t,0})\|_2^2 \to 0,
\end{align*}
as $t \to \infty$. Note that 
$$
\|\nabla \J(z^{t},w^{t})\|_2^2 = \|\nabla_z \J(z^{t},w^{t})\|_2^2 + \|\nabla_w \J(z^{t},w^{t})\|_2^2,
$$
and, thus, it remains to show that $\|\nabla_w \J(z^{t},w^{t})\|_2^2 \to 0$ for $t \to \infty$. Using triangle inequality, we obtain
$$
0 \le \|\nabla_w \J(z^{t},w^{t}) \|_2 = \|\nabla_w \J(z^{t},w^{t}) - \nabla_w \J(z^{t},w^{t-1})\|_2 + \| \nabla_w \J(z^{t},w^{t-1}) \|_2
$$
We already showed that the second summand converges to zero. For the first summand we can use the fact that $\nabla_w \J(z,w)$ is continuous for every $\varepsilon>0$. Therefore, the first summand converges to $0$ if $\| w^t - w^{t-1}\|_2 \to 0$ for $t\to \infty$. In fact, we have
\begin{align*}
0 \le \| w^t - w^{t-1}\|_2 & = \| w^{t,I_w} - w^{t,0}\|_2  
=  \big\| \sum_{j = 1}^{I_w} \nu_{t,j} \nabla_w \J(z^{t},w^{t,j-1}) \big\|_2 \\
& \le \sum_{j = 1}^{I_w} \nu_{t,j} \| \nabla_w \J(z^{t},w^{t,j-1}) \|_2 \\
& \le \tau^{-N} \beta_T^{-1} \sum_{j = 1}^{I_w} \| \nabla_w \J(z^{t},w^{t,j-1}) \|_2,
\end{align*}
where in the last line we used Lemma \ref{l: learning rates finite}. 
Taking $t \to \infty$, we obtain \mbox{$\| w^t - w^{t-1}\|_2 \to 0$} and, consequently, by continuity 
\[
\|\nabla_w \J(z^{t},w^{t}) - \nabla_w \J(z^{t},w^{t-1})\|_2 \to 0,\] 
which gives 
$
\|\nabla_w \J(z^{t},w^{t}) \|_2 \to 0,
$
as $t \to \infty$. 

For the convergence speed, we consider the sequence
$$
s_t := \max\big\{ \min_{i =1 ,\dots, I_z} \|\nabla_z \J(z^{t,i-1},w^{t-1})\|_2^2 , \min_{j =1 ,\dots, I_w} \|\nabla_w \J(z^{t},w^{t,j-1})\|_2^2    \big\}.
$$
If $s_t$ is small, it implies that the gradient iterations in either direction are small and we are in proximity of stationary point. For the minimum of $s_t$ the following upper bound holds,
\begin{align*}
& \min_{t = 1, \dots, T} s_t 
\le \tfrac{1}{T} \sum_{t= 1}^{T} s_t  \\
& \quad \quad \le \tfrac{1}{T} \sum_{t= 1}^{T} \big[ \min_{i =1 ,\dots, I_z} \|\nabla_z \J(z^{t,i-1},w^{t-1})\|_2^2 + \min_{j =1 ,\dots, I_w} \|\nabla_w \J(z^{t},w^{t,j-1})\|_2^2 \big].
\end{align*}

The first minimum is bounded from above as
\begin{align*}
\min_{i =1 ,\dots, I_z} \|\nabla_z \J(z^{t,i-1},w^{t-1})\|_2^2
& \le \frac{1}{I_z} \sum_{i=1}^{I_z} \|\nabla_z \J(z^{t,i-1},w^{t-1})\|_2^2 \\
& \le \frac{1}{I_z \ds\min_{i = 1 \dots, I_z} \mu_{t,i} } \sum_{i=1}^{I_z} \mu_{t,i} \|\nabla_z \J(z^{t,i-1},w^{t-1})\|_2^2.
\end{align*}
and, similarly, the second minimum is bounded by
$$
\min_{j =1 ,\dots, I_w} \|\nabla_w \J(z^{t},w^{t,j-1})\|_2^2
\le \frac{1}{I_w \ds\min_{i = 1 \dots, I_w} \nu_{t,j} } \sum_{j=1}^{I_w} \nu_{t,j} \|\nabla_w \J(z^{t},w^{t,j-1})\|_2^2.
$$
Combined with \eqref{eq: objective decreas}, these bounds give us
\begin{align*}
\min_{t = 1, \dots, T} s_t 
& \le \frac{1}{T \min\{I_z \ds\min_{i = 1 \dots, I_z} \mu_{t,i} , I_w \ds\min_{j = 1 \dots, I_w} \nu_{t,j} \} }  \J(z^0, w^0) \\
& \le \frac{1}{T \min\{I_z, I_w \} \min \{ \ds\min_{i = 1 \dots, I_z} \mu_{t,i} ,\ds\min_{j = 1 \dots, I_w} \nu_{t,j}  \} }  \J(z^0, w^0).
\end{align*}
In a view of \eqref{eq: lower bound learning rate} and \eqref{eq: upper bound window}, for the step sizes $\mu_{t,i}$ we obtain
\begin{align*}
\frac{1}{\ds\min_{i = 1 \dots, I_z}{\mu_{t,i}} } & 
\le \frac{1}{\mu_{t,c}} \le B_w = \beta_T^{-1} \J(z^0,w^0) \max_{\ell = 1, \dots, L} \lambda_\ell^{-2} \| F_\ell \|^2  + \alpha \\ 
& \le \max\{\alpha_T^{-1}, \beta_T^{-1} \} \J(z^0,w^0) \max_{\ell = 1, \dots, L} \lambda_\ell^{-2} \| F_\ell \|^2 + \max\{ \alpha, \beta \} ,
\end{align*}
and for $1/\ds\min_{j = 1 \dots, I_w}{\nu_{t,j}}$ the bounds is precisely the same. Then, we arrive at
\begin{align*}
\min_{t =1 ,\dots, T} s_t
& \le \frac{\max\{ \alpha_T^{-1}, \beta_T^{-1} \} (\J(z^0,w^0))^2 \ds\max_{\ell = 1, \dots, L} \lambda_\ell^{-2} \| F_\ell \|^2 + 
\J(z^0,w^0) \max\{ \alpha, \beta \}
}{T \min\{I_z, I_w \} }.
\end{align*}
\bend

We note that Algorithm \ref{alg:PP1} will always converge in terms of the value of the objective function, even if the condition $\alpha_T, \beta_T > 0$ is violated. That is the sequence $( \J(z^t, w^t) )_{t \ge 0}$ is a bounded from below and non-increasing sequence. Hence, there exits limit $\lim_{t \to \infty} \J(z^t, w^t)$. Moreover, the continuity of $\J$ implies the existence of limit points $\lim_{t \to \infty} (z^t,w^t) = (z_*, w_*)$. However, we cannot guarantee that the gradient vanishes at the limit points.
The condition $\alpha_T, \beta_T > 0$ is also related to the fact that the loss function $L_\varepsilon$ is invariant to the rescaling of the object and the mask, so that for all $a \in \C$ equality $L_\varepsilon(z;w) = L_\varepsilon(a z;w/a)$ holds, while the gradients $\nabla_z L_\varepsilon, \nabla_w L_\varepsilon$ are scale-dependant. The condition $\alpha_T, \beta_T > 0$ includes Tikhonov regularization terms to the objective $\J$, which partially resolves this ambiguity to cases when $|a| = 1$, a so-called global phase factor.

\section{Numerical Examples}\label{sec:NuEx}

\subsection{Experimental setup}
\begin{figure}[b!]
\centering
\includegraphics[width=1\textwidth]{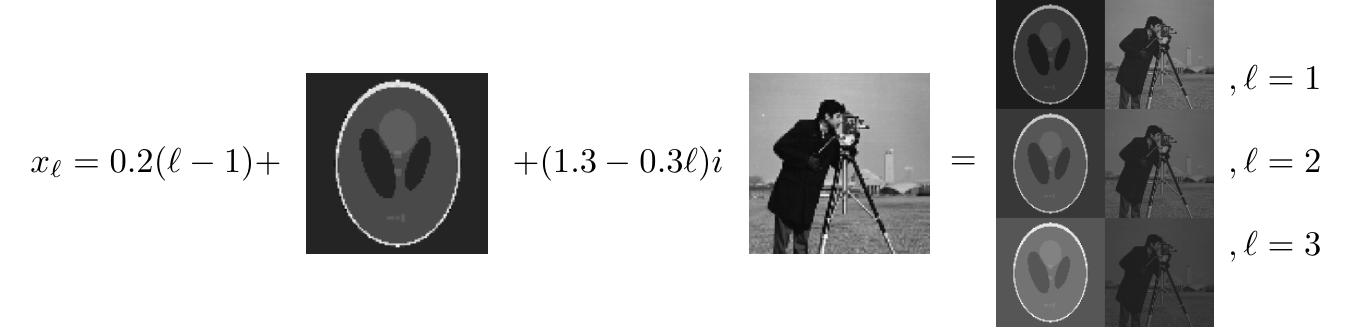}
\caption{Synthetic object in polychromatic light, $d = 100\times100$, $L=3$.}\label{fig:3}
\end{figure}

In this section we perform numerical experiments to explore the performance of gradient descent for polychromatic ptychography. 
All our experiments will be performed on a synthetic data within the following setup. We will consider the polychromatic light with $L=3$ wavelengths $\lambda = (1,1.25,1.5)$. 
For the object $z$ an image of size $d = 100\times100$ is used, where the real and the imaginary parts are scaled images of the Shepp-Logan phantom and the cameraman, respectively. We slightly alternate the real and imaginary part for different wavelengths $\lambda$ to imitate the dependence of the object on the wavelength as described in Figure \ref{fig:3}.     

\begin{figure}
\centering
\includegraphics[width=0.5\textwidth]{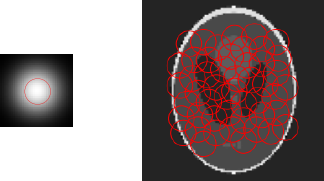}
\caption{Localized mask of size $40\times40$ and its shifts. Each of 49 red circles indicates the position of the mask on the Fermat spiral \eqref{eq: Fermat spiral}.}\label{fig:4}
\end{figure}

The mask $w$ is assumed to be locally supported with $\supp(w) = [\delta]^2$, $\delta = 40$. Within the support, its values are sampled from the Gaussian function,
$$
g_k = e^{ - \| k  - \mu \|_2^2 / 2 \sigma^2 }, \quad k \in [\delta]^2
$$
with $\mu = ((\delta-1)/2, (\delta-1)/2)$ and $\sigma^2 = \delta^2/ 20$. Then  the window is formed by rescaling $w_\ell = \sqrt{\sigma_\ell} g / \| g \|_2$, where weights $\sigma = (0.2,0.5,0.3)$ represent the spectral density of polychromatic light. For the visualization, we refer reader to Figure \ref{fig:4}.

The set of shifts is selected by moving the center of the mask along the Fermat spiral as discussed in \cite{Huang.2014}. That is, in polar coordinate system $(r,\phi)$ the center of the mask satisfies
\begin{equation}\label{eq: Fermat spiral}
r_k = c_{sp} \sqrt k, \ \phi_k = k \phi_0, \ 0 \le k \le \lceil 0.5((d - \delta)/ c_{sp})^2 \rceil 
\end{equation}
where the $c_{sp} = 4.9$ is the scaling factor of the radius and the initial angle $\phi_0$ is given by $\phi_0 = 2\pi(\frac{2}{1 + \sqrt 5})^2 \approx 137.508^\circ$. Then, pairs $(r_k,\phi_k)$ are transformed into Cartesian coordinate system as $m = (m^1_k,m^2_k)$ and $\mathcal M$ only contains those points $m = (m^1_k + d/2,m^2_k + d/2)$, for which non-zero entries of $w$ are contained inside the object domain $[d]^2$ as depicted in Figure \ref{fig:4}.

The measurements $y_{m,k}$, resulting from simulated polychromatic ptychographic experiment are represented by the weighted sum of intensities as in equation \eqref{eq: poly meas} and Figure \ref{fig:2}.
Furthermore, the measurements are corrupted by the Poisson noise, so that
$$
\tilde y_{m,k} = \frac{d}{N_{p}} \text{Pois}\left(\frac{N_{p} y_{m,k}}{d}\right),
$$
where $N_{p} = 10^{6}$ represents the number of photons used for the experiment. Since we fix the random seed for reproducible results, the relative noise level $\| \sqrt{ \tilde y} - \sqrt{y}\|_2/ \| \sqrt{y} \|_2$ for all experiments is $\approx 0.1024$. 

As a measures of performance we will consider loss functions $L_\varepsilon$ defined in equation \eqref{eq:Re1} with $\varepsilon = 10^{-8}$. In addition, with true object $z$ known for synthetic data, the relative object error 
\begin{equation}
\| z - z^t \|_2 / \| z \|_2  \label{eq: total relative error}\\
\end{equation}
can be evaluated and will be used for comparisons.

The proposed algorithms are compared to ptychographical information multiplexing method (PIM$_\alpha$) with parameter $\alpha$ denoting the step size \cite{Batey.2014}. 

We note that all experiments were performed in Python on the laptop running Windows 10 Pro with an Intel(R) Core(TM) i7-8550U processor and with 16 GB RAM. 

\subsection{Non-blind polychromatic ptychography
}
We start with the non-blind polychromatic ptychography. In order to reconstruct the object, we perform the gradient descent minimization of the loss function $J$ with different parameters. For the first reconstruction denoted by AF$_{0}$, the step size is constant $\mu_t = \mu_c = 1/B(w)$, where $B(w)$ as in Corollary \ref{col: non-blind ptycho} and regularization parameters are set $\alpha_T = 10^{-2}$ and $\alpha_S =0$. In second trial we additionally include smoothness penalty $\alpha_S = 0.1$ to highlight its benefit. This algorithm is denoted by AF$_{0.1}$. In addition to smoothness penalty, the third recovery procedure AF$_{0.1}$+AGA selects the step size via AGA with $N=1$ and $\tau = 0.5$. Last algorithm is PIM$_\alpha$ with $\alpha = 1$. Since the window is known, we only use object update of PIM. For all algorithms an initial guess $z^0$ is the flat object, that is $(z^0)_k = 1$ for all $k \in [d]^2$. 

\begin{figure}[b!]
\centering
\includegraphics[width=1\textwidth]{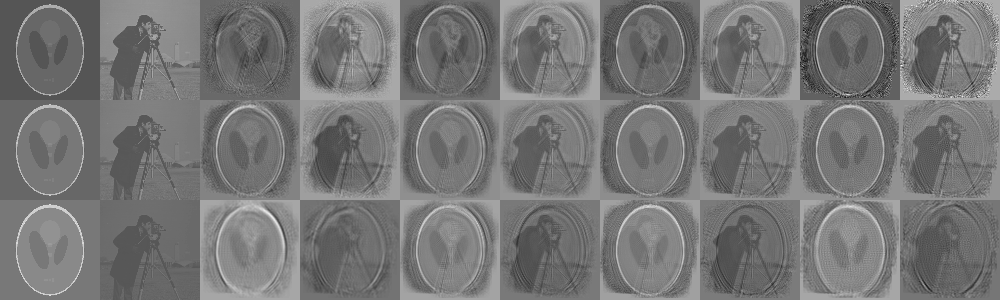}
\caption{Reconstruction of the object with the known mask. Each row corresponds to a single wavelength $\ell =1,2,3$. The two consecutive columns are the real and imaginary parts of the object. In figure, we show the true object $z$ and reconstructed objects with AF$_0$, AF$_{0.1}$, AF$_{0.1}$+AGA and PIM$_1$. 
}\label{fig:5}
\end{figure}

\begin{figure}[t!]
\centering
\begin{minipage}{0.45\textwidth}
\subfloat[Loss function $L_\varepsilon$. \label{fig:61}]
{\includegraphics[width = 1\linewidth]{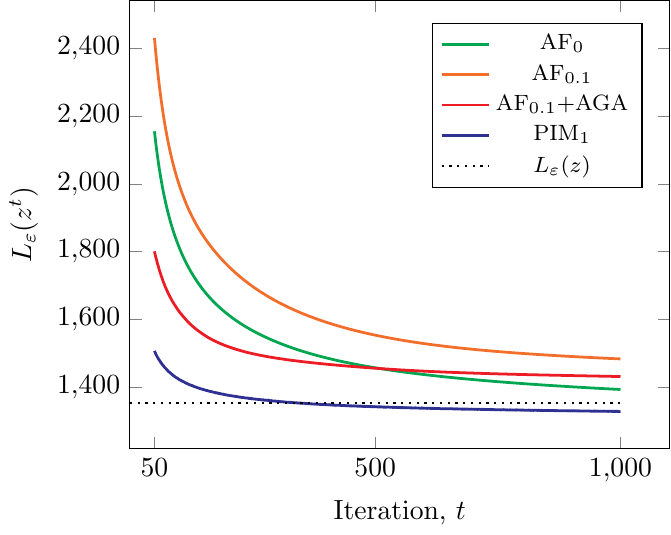}}
\end{minipage}
\hspace{2mm}
\begin{minipage}{0.45\textwidth}
\subfloat[Relative object error.\label{fig:62}]
{\includegraphics[width = 1\linewidth]{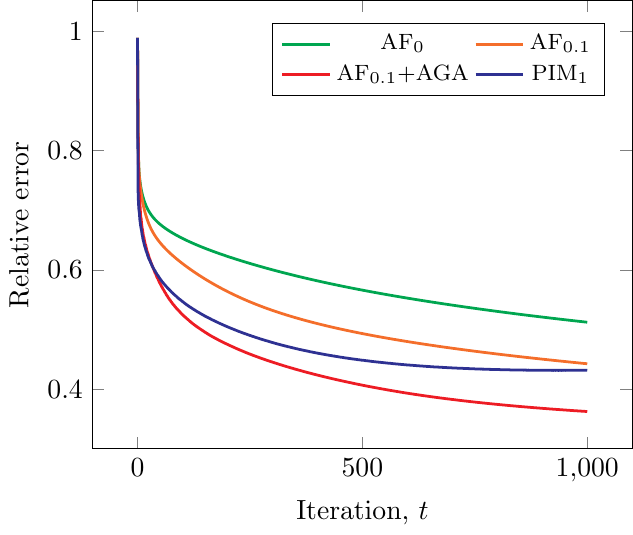}}
\end{minipage}
\caption{Numerical comparison of reconstructions with different algorithms. First $50$ iterations are excluded for better visualization.
}\label{fig:6}
\end{figure}

The outcome of the $1000$ iterations of each algorithms are presented in Figure \ref{fig:5}. Furthermore, the numerical comparison between the methods is shown in Figure \ref{fig:6}. 

We observe the difference between reconstructions with AF$_0$ and AF$_{0.1}$, which suggests that for continuous in $\lambda$ objects the use of smoothness penalty is beneficial. Furthermore, inclusion of AGA significantly speeds up the convergence of gradient descent. Comparing the reconstructions with two non-regularized algorithms, AF$_0$ and PIM$_1$, the latter produces visually better result and is the fastest in minimizing $L_\varepsilon$. However, both non-regularized methods provide a smaller final value of the loss $L_\varepsilon$, the relative error is larger compared to AF$_{0.1}$+AGA. While the difference is striking for AF$_0$ and AF$_{0.1}$, for PIM$_1$ high frequency artifacts can be observed, especially prominent for the third wavelength.  This may hint towards an occurrence of the overfitting phenomena, which may be prevented by an inclusion of the smoothness penalty.   

The runtime of AF$_0$ and AF$_{0.1}$ is 6m 53s and 6m 21s, respectively. The use of AGA leads to additional evaluations of the objective function and slightly longer runtime of 8m 2s. For PIM it took 8m 40s to perform 1000 iterations. 

\subsection{Blind polychromatic ptychography}
In the next experiment, we assume that the mask $w$ is unknown. Then, a reconstruction of both the object and the window is performed via Algorithm \ref{alg:PP1} and two versions of PIM$_\alpha$, with $\alpha=0.1$ and $\alpha = 1$.
For Algorithm \ref{alg:PP1}, denoted by $AF$, the number of iterations is set to $T = 1000$ with the object and mask iterations $I_z = I_w = 1$. This corresponds to $1000$ gradient steps for each the object and the mask. The object regularization parameters $\alpha_T = 10^{-2}$ and  $\alpha_S=0.1$ are set as for non-blind experiment above and mask regularization parameters are set to $\beta_T = 10$ and $\beta_S = 0$. The step size is selected via AGA with $N=1$ and $\tau = 0.5$. We note that $\beta_S$ is set to zero as components of the window have different norms corresponding to the spectral density of the light distribution. Thus, it is not expected that the window components $w_\ell$ should be close to each other. 

For the object initialization the flat starting point $z^0$ is used and the initial guess for the mask is given by $w_\ell^0 = L^{-1/2} g^0 / \|g^0 \|_2$ with
$$
g^0_{k} =
\begin{cases}
2.3, & \| k - \mu \|_2 \le \sqrt{0.3} \delta/2,\\
1.3, & \| k - \mu \|_2 \le \sqrt{0.6} \delta/2,\\
0.3, & \| k - \mu \|_2 > \sqrt{0.6} \delta/2, k \in [\delta]^2,\\
0, & \text{otherwise}.
\end{cases}
$$

\begin{figure}[b!]
\centering
\includegraphics[width=1\textwidth]{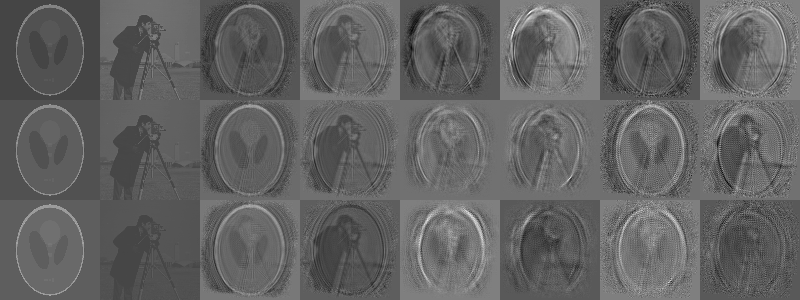}
\caption{Reconstruction of the object for blind polychromatic ptychography. Each row corresponds to the single wavelength $\ell =1,2,3$. The two consecutive columns are the real and imaginary parts of the object. In figure, we the show true object $z$ and reconstructed objects with AF, PIM$_{0.1}$ and PIM$_1$. 
}\label{fig:8}
\end{figure}

The motivation behind this construction is to provide the rough approximation of the shape of the true mask $w$ with the energy equally distributed along the wavelength spectrum. As the performance of the gradient methods applied to a non-convex functions are known to be sensitive to the initial guess, this initialization is hopefully sufficient to ensure the fast convergence to the true mask. 

\begin{figure}[t!]
\centering
\includegraphics[width=1\textwidth]{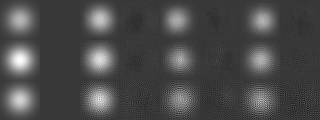}
\caption{Reconstruction of the window for blind polychromatic ptychography. Each row corresponds to the single wavelength $\ell =1,2,3$. The two consecutive columns are the real and imaginary parts of the mask. The two consecutive columns are the real and imaginary parts of the window. In figure, we the show true window $w$ and reconstructed objects with AF, PIM$_{0.1}$ and PIM$_1$.  
}\label{fig:9}
\end{figure}

\begin{figure}[b!]
\centering
\subfloat[Loss $L_\varepsilon$. \label{fig:101}]
{\includegraphics[width = 0.315\textwidth
]{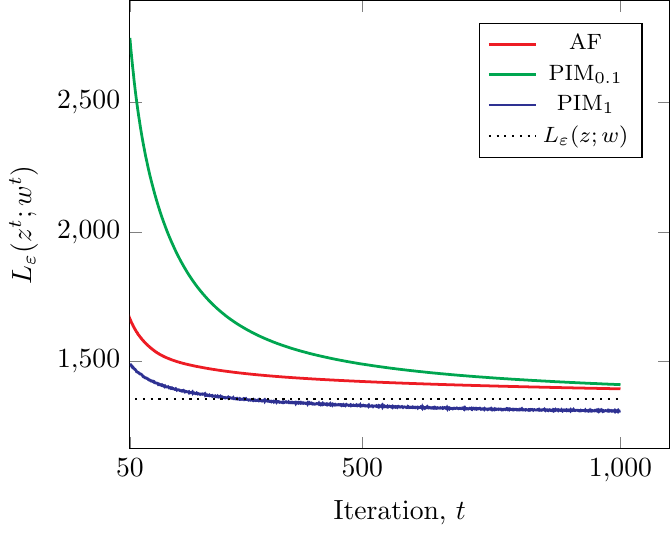}}
\hspace{2mm}
\subfloat[Object error. \label{fig:102}]
{\includegraphics[width = 0.30\textwidth]{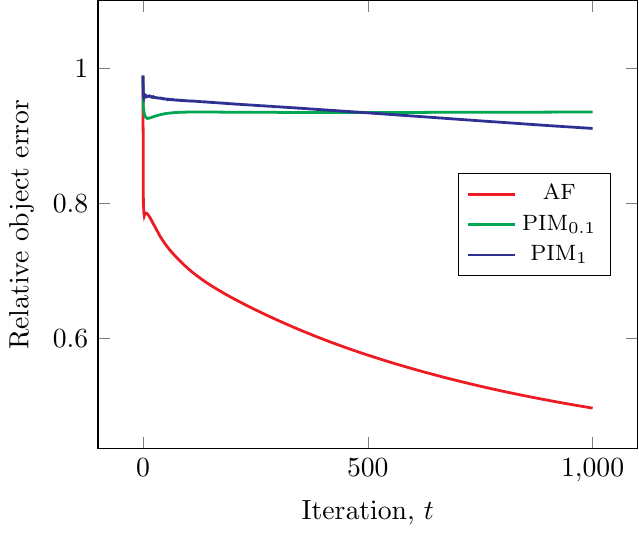}}
\hspace{2mm}
\subfloat[Window error. \label{fig:103}]
{\includegraphics[width = 0.30\textwidth]{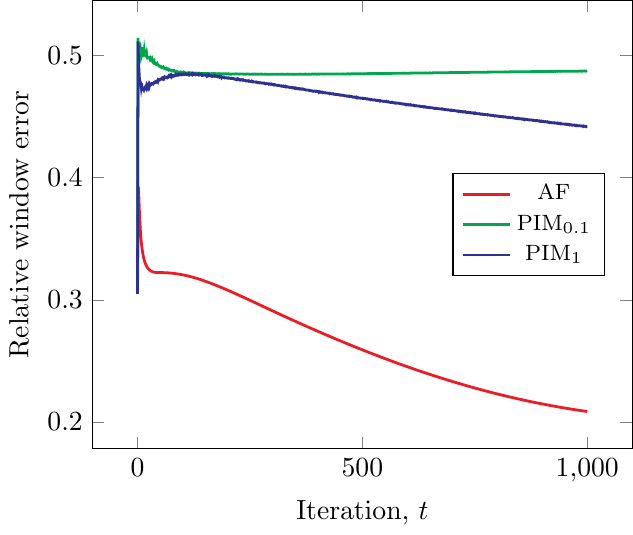}}
\caption{Numerical comparison of AF and PIM for blind polychromatic ptychography. First $50$ iterations are excluded for better visualization.
}\label{fig: 10}
\end{figure}

Reconstructions with selected algorithms are shown in Figures \ref{fig:8} and \ref{fig:9}. We observe that reconstruction with AF is visually closer to the original image. Figure \ref{fig:101} shows that the values of the loss function $L_\varepsilon$ for PIM$_1$ are smaller than for AF, which point towards the overfitting again. 

It it also notable that the reconstruction for the last wavelength is the most noisy. Recall that the measurements are a mixture of the intensities for each wavelength weighted with $\lambda_\ell^{-2}$. Furthermore, each intensity is scaled with the spectral density distribution $\sigma_\ell$ of the corresponding wavelength contained in $w_\ell$. Therefore, the weights  $\sigma_\ell/\lambda_\ell^{2} = (0.2,0.32,0.1333)$ imply that the reconstruction for the last wavelength is, roughly, three times more sensitive to the noise then the second.

The runtime of AF is 20m 50s, while PIM only requires 7m 52s and 7m 17s. Compared to non-blind case AF performs double the number of gradient steps, which why the runtime was expected to at least double from 6m 21s to $\approx$13m. The extra 8m result from the recomputation of the step sizes for object and window in Algorithm \ref{alg:PP1}. 


\begin{figure}[t!]
\centering
\includegraphics[width=1\textwidth]{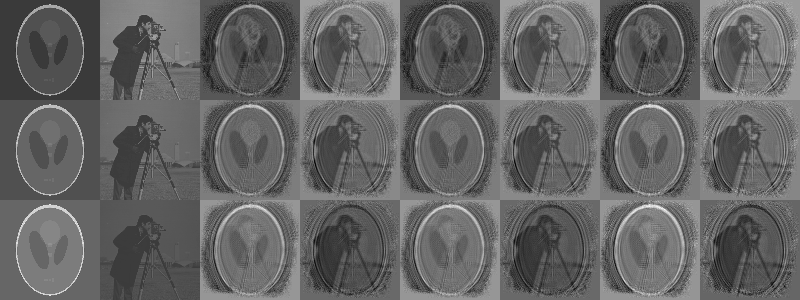}
\caption{Reconstruction of the object for blind polychromatic ptychography with different iteration number setups. Each row corresponds to the single wavelength $\ell =1,2,3$. The two consecutive columns are the real and imaginary parts of the mask. In figure, we show the true object $z$ and reconstruction with Algorithm \ref{alg:PP1} for three choices of $(T,I_z,I_w) = (200,10,10), (400,5,5)$ and  $(2000,1,1)$. 
}\label{fig:11}
\end{figure}

\begin{figure}[b!]
\centering
\subfloat[Loss $L_\varepsilon$. \label{fig:121}]
{\includegraphics[width = 0.315\textwidth
]{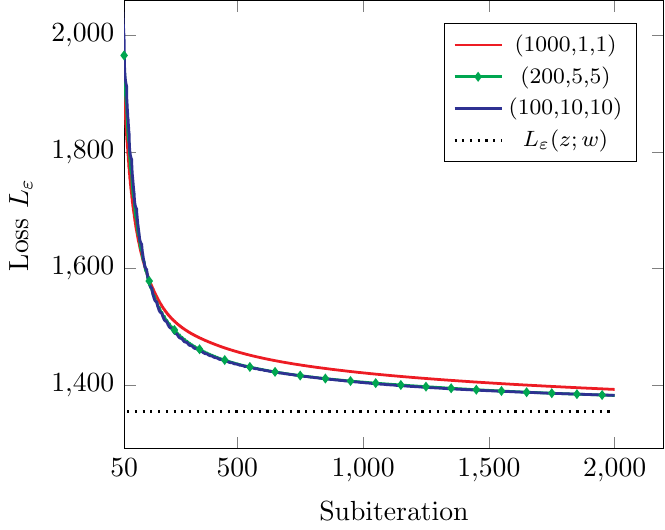}}
\hspace{2mm}
\subfloat[Object error. \label{fig:122}]
{\includegraphics[width = 0.30\textwidth]{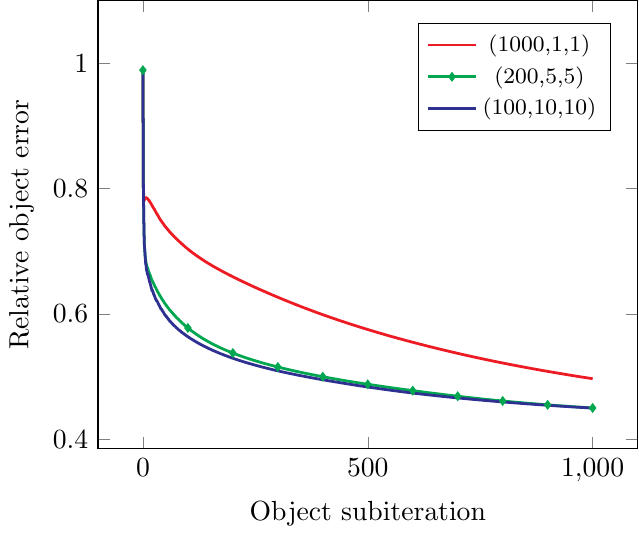}}
\hspace{2mm}
\subfloat[Window error. \label{fig:123}]
{\includegraphics[width = 0.30\textwidth]{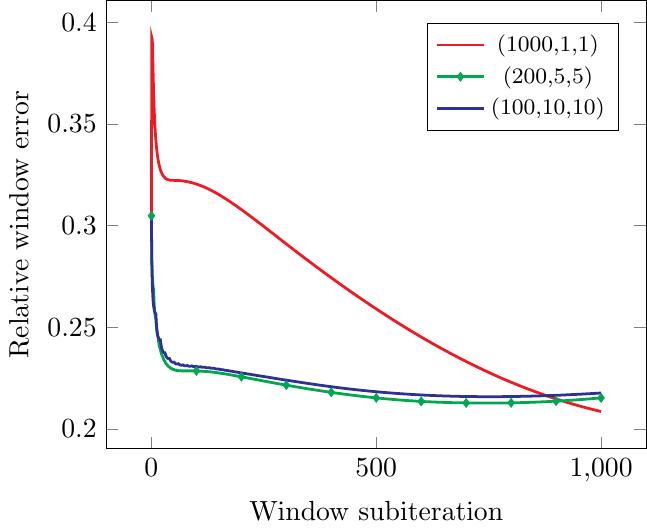}}
\caption{Numerical comparison of Algorithm \ref{alg:PP1} for three choices of $(T,I_z,I_w) = (200,10,10), (400,5,5) , (2000,1,1)$. First $50$ iterations are excluded for better visualization.
}\label{fig: 12}
\end{figure}

Therefore, we explore if increasing the number of object and window subiterations $I_z$ and $I_w$ will improve the runtime without interfering the quality of reconstruction. Thus, we repeat the reconstruction for $(T,I_z,I_w) = (200,5,5)$ and $(T,I_z,I_w) = (100,10,10)$. For all three sets of parameters the number of gradient steps for both the object and the mask remains $T I_z = T I_w = 1000$. As the Figure \ref{fig:11} shows, the redistribution of the iterations have visual impact on the reconstruction quality. In Figure \ref{fig: 12}, we observe a larger number of iterations leads to a faster decay of the loss function and the errors. Note that in Figure \ref{fig:123} relative window error for set-up $(1000,1,1)$ increases as after first object subiteration, the window gradient points into wrong direction while for $(200,5,5)$ and $(100,10,10)$ larger number of object subiterations leads to a better direction for the window. The runtimes of the algorithms are 20m 50s, 14m 42s, 14m 37s for parameters $(T,I_z,I_w) = (1000,1,1), (200,5,5), (100,10,10)$, respectively.

\section{Conclusion and Discussion}

In this paper, we consider recovery from polychromatic ptychographic measurements. It is performed via gradient descent applied to the amplitude-based squared loss with guaranteed convergence to a critical point. We also combine the idea of gradient-based optimization with alternating minimization to address blind polychromatic ptychography. This results in Algorithm \ref{alg:PP1} with convergence guarantees to a critical point summarized in Theorem \ref{theo:PP1}. 

We note that polychromatic ptychography can be seen as a generalization of single wavelength ptychography as two measurement models coincide if $L$ is set to one. Therefore, our analysis generalized some results in the literature. In particular, for the non-blind ptychography gradient descent for the loss function $L_\varepsilon$ was already studied in \cite{Xu.2018}. Furthermore, alternating minimization is sometimes used for blind ptychography \cite{Hesse.2015,Chang.2019,Fannjiang.2020} and these works could be compared to our results. However, to our knowledge none of the methods for blind ptychography guarantees a sublinear convergence rate as in Theorem \ref{theo:PP1}. The proofs presented in this paper can be extended for layerwise optimization algorithm for multislice ptychography \cite{Bangun.2022}. 

While Algorithm \ref{alg:PP1} is supported by theoretical analysis, numerical examples point towards its underperformance in terms of computation time. This problem can be tackled by the use of high performance computing, e.g. parallelized computation of gradients.
Another potential way to computational efficiency is transition to alternating stochastic gradient descent, convergence of which would combine proof ideas of Theorem \ref{theo:PP1} and \cite{Melnyk.2022}. 

In numerical trials, we also observed that reconstruction is sometimes contains high frequency noise. Such artifacts could be avoided by inclusion of additional regularized such as smoothness or total variation penalties.  

\section*{Acknowledgments}
This work was funded by the Helmholtz Association under contracts No.~ZT-I-0025 (Ptychography 4.0), No.~ZT-I-PF-4-018 (AsoftXm), No.~ZT-I-PF-5-28 (EDARTI).

\bibliography{cit}  

\end{document}